\documentclass[10pt]{amsart}

\usepackage{amsmath,amssymb,amsthm, amsfonts}
\usepackage{array,mathtools} 
\usepackage{mathrsfs}
\usepackage{euler}
\usepackage{times}

\usepackage[color,matrix,arrow,curve]{xy}
\usepackage{tikz-cd}
\usepackage{hyperref}

\theoremstyle{plain}

\newtheorem{thm}{Theorem}[section]
\newtheorem{prop}[thm]{Proposition}
\newtheorem{lem}[thm]{Lemma}
\newtheorem{lemma}[thm]{Lemma}
\newtheorem{cor}[thm]{Corollary}
\newtheorem{set}[thm]{Setting}
\newtheorem{ass}[thm]{Assumption}

\theoremstyle{remark}
\newtheorem{rmk}[thm]{Remark}

\newcommand{\bC}{\mathbb{C}}

\newcommand{\bF}{\mathbb{F}}
\newcommand{\bG}{\mathbb{G}}
\newcommand{\bL}{\mathbb{L}}

\newcommand{\bP}{\mathbb{P}}
\newcommand{\bQ}{\mathbb{Q}}
\newcommand{\bZ}{\mathbb{Z}}

\newcommand{\cE}{\mathcal{E}}
\newcommand{\cF}{\mathcal{F}}
\newcommand{\cG}{\mathcal{G}}
\newcommand{\cO}{\mathcal{O}}
\newcommand{\cI}{\mathcal{I}}
\newcommand{\cM}{\mathscr{M}}
\newcommand{\cN}{\mathcal{N}}

\newcommand{\sE}{\mathscr{E}}
\newcommand{\sM}{\mathscr{M}}

\newcommand{\qcoh}{\mathrm{D}_\mathrm{qcoh}}
\newcommand{\coh}{\mathrm{D}_\mathrm{coh}}
\newcommand{\perf}{\mathrm{D}_\mathrm{perf}}

\DeclareMathOperator{\rhom}{\mathbf{R}\mathcal{H}om}

\newcommand{\at}{\mathrm{At}}
\newcommand{\At}{\mathrm{at}}
\newcommand{\spec}{\mathrm{Spec}}
\newcommand{\rmd}{\mathrm{D}}
\newcommand{\pr}{\mathrm{pr}}
\newcommand{\pt}{\mathrm{pt}}

\newcommand{\tensorhom}{\otimes\text{-}\mathrm{hom}}
\newcommand{\homtensor}{\mathrm{hom}\text{-}\otimes}
\newcommand{\tensor}{\otimes\text{-}\otimes}

\newcommand{\lst}[1]{#1^*\!\text{-}\!\ #1_*}
\newcommand{\usk}[1]{#1_*\!\text{-}\!\ #1^!}

\newcommand{\otv}[1]{[#1]^\mathrm{vir}_\mathrm{OT}}
\newcommand{\bfv}[1]{[#1]^\mathrm{vir}_\mathrm{BF}}

\newcommand{\Supp}{\mathop{\rm Supp}\nolimits}
\newcommand{\Pic}{\mathop{\rm Pic}\nolimits}
\newcommand{\id}{\mathrm{id}}
\newcommand{\ch}{\mathop{\rm ch}\nolimits}
\newcommand{\td}{\mathop{\rm td}\nolimits}
\newcommand{\ev}{\mathop{\rm ev}\nolimits}
\newcommand{\HOM}{\mathop{\mathcal{H}om}\nolimits}
\newcommand{\tr}{\mathop{\rm tr}\nolimits}

\newcommand{\DT}{\mathop{\rm DT}\nolimits}

\title[GV invariants of fiber classes on CY 4-folds fibered over curves]{Gopakumar-Vafa invariants of fiber classes on Calabi-Yau 4-folds fibered over curves}
\date{}

\author{Yalong Cao}
\address{Morningside Center of Mathematics, Institute of Mathematics \& State Key Laboratory of Mathematical Sciences,
Academy of Mathematics and Systems Sciences,
 Chinese Academy of Sciences, 
 Beijing, 
 China}
\email{yalongcao@amss.ac.cn}

\author{Feng Qu}
\address{Northeast Normal University, School of Mathematics and Statistics, Changchun, Jilin, China}
\email{quf996@nenu.edu.cn}

\begin{document}

\begin{abstract}
We prove a conjectural correspondence of Cao-Maulik-Toda which 
relates Gopakumar-Vafa invariants of fiber classes on a smooth projective Calabi-Yau 4-fold fibered over a curve to the Gopakumar-Vafa invariants of a smooth fiber 
under an orientation compatibility assumption on the moduli spaces.
\end{abstract}

%\keywords{}
%\subjclass{}
% ${}$ \\ \textbf{Keywords}: Tautological bundles, Hilbert scheme of points, Calabi-Yau 4-folds, virtual pullback
% ${}$ \\ \textbf{MSC 2010}: 14N35, 14J32, 14C17

\maketitle

%\tableofcontents

\section{Introduction}

\subsection{Gopakumar-Vafa  invariants on Calabi-Yau 4-folds}

Let $X$ be a smooth projective Calabi-Yau 4-fold over $\mathbb{C}$ and $\beta\in H_2(X,\mathbb{Z})$ a non-zero curve class. 
Klemm and Pandharipande \cite{KP} introduced \textit{genus zero and genus one BPS invariants}: 
\begin{equation}\label{intro equ on gv}
n_{0,\beta}(\gamma)(X),  \quad n_{1,\beta}(X), 
\end{equation}
defined from Gromov-Witten theory, where $\gamma\in H^4(X,\mathbb{Z})$. From virtual dimension counting, Gromov-Witten invariants on projective Calabi-Yau 4-folds
are non-trivial only when genus $g\leqslant 1$, so are the BPS invariants.
The invariants \eqref{intro equ on gv} are conjectured to be integers and are analogues of Gopakumar-Vafa invariants on Calabi-Yau 3-folds \cite{GV}, 
and we will refer to them as \textit{Gopakumar-Vafa (type) invariants of Calabi-Yau 4-folds}. 

In \cite{CMT1, CT2}, the authors gave a sheaf theoretic interpretation of \eqref{intro equ on gv} using the moduli scheme $M_{\beta}(X)$ of one dimensional
stable sheaves $\cF$ with $[\cF]=\beta$ and $\chi(\cF)=1$. There is a $\DT_4$ virtual class \cite{BJ} (see also \cite{CL1} in special cases)
\[
[M_{\beta}(X)]^{\mathrm{vir}} \in H_2(M_{\beta}(X),\mathbb{Z}), 
\]
which depends on the choice of orientation \cite{CGJ} on the moduli space. On each connected component, there are two choices of orientation. Reversing the choice 
affects the virtual class by a minus sign componentwise. 
The recent work of Oh-Thomas \cite{OT} lifts the $\DT_4$ virtual class from
homology to Chow groups after inverting two in coefficients. 

We define cohomology classes 
$$
\tau_k\colon H^{m}(X,\bQ)  \to  H^{2k+m-2}(M_\beta(X),\bQ) 
$$  
by Chern characters  of the `normalized' universal sheaf and Fourier-Mukai type construction \eqref{equ on tau}.
The corresponding $\DT_4$ invariants are 
\begin{equation}\label{intro equ on tau inv}
\langle \tau_k(\gamma) \rangle_{X,\beta} := \int_{ [M_{\beta}(X)]^{\mathrm{vir}} } \tau_k(\gamma) \in \bQ. \end{equation}
Sometimes we omit the subscript $X$  if there is no confusion. 

Conjecturally, \eqref{intro equ on tau inv} recovers \eqref{intro equ on gv} \cite{CMT1, CT2}: e.g., in the genus zero case, 
it is conjectured that 
$$
\langle \tau_0(\gamma)\rangle_{X,\beta}=n_{0,\beta}(\gamma)(X), 
$$ 
for some choice of orientation on the LHS.
We refer to \cite{CK2, CKM1, CKM2, CMT2, CT1,CT3,CT4} for related studies on Gopakumar-Vafa invariants \eqref{intro equ on gv}, \cite{COT1,COT2} for 
discussions of the holomorphic symplectic case, and \cite{CKM3} for extensions to the orbifold case.

\subsection{Calabi-Yau 4-folds fibered over curves}

Let $X$ be a connected smooth complex projective CY 4-fold and 
$$\pi \colon X \to C$$ 
be a projective surjective map to a smooth projective curve $C$ with connected fibers.

By generic smoothness of $\pi$, there exists a non-empty open subset $U\subset C$ such that 
$$\pi: \pi^{-1}(U)\to U$$
is smooth. For $c\in U$, the fiber $X_c:=\pi^{-1}(c)$ is a smooth projective Calabi-Yau 3-fold. If $\beta'\in H_2(X_c,\bZ)$ is a non-zero curve class, then
one has \textit{genus zero Gopakumar-Vafa invariant}
$$n_{0,\beta'}(X_c). $$
There is also a sheaf theoretic interpretation of this invariant using the moduli scheme $M_{\beta'}(X_c)$ of one dimensional
stable sheaves $\cF$ on $X_c$ with $[\cF]=\beta'$ and $\chi(\cF)=1$. 
It is conjectured in \cite{HST, Katz} that
\begin{equation}\label{equ on kz} n_{0,\beta'}(X_c)=\int_{[M_{\beta'}(X_c)]^{\mathrm{vir}}}1, \end{equation}
where $[M_{\beta'}(X_c)]^{\mathrm{vir}}$ is the $\DT_3$ virtual class \cite{BF1,Tho}. 

Let $\beta\in H_2(X,\mathbb{Z})$ be a fiber class,~i.e., $\pi_*\beta=0$, and $M_\beta(\pi)$ denote the moduli scheme of one-dimensional stable sheaves $\cF$ on the fibers of $\pi$ with $[\cF]=\beta$ and $\chi(\cF)=1$. By Proposition~\ref{relabs}, there is a canonical isomorphism 
$$M_\beta(\pi) \simeq M_\beta(X). $$
Let $M_\beta(X_c)$ be the fiber of $M_\beta(\pi) \to C$ over $c$, which fits into the Cartesian diagram
\begin{align*}
\xymatrix{
M_\beta(X_c) \ar[r]^{ } \ar[d]  \ar@{}[dr]%|{\Box} 
	& M_\beta(\pi) \ar[d] \\
c \ar[r]^{i_c} 
	&  C.}
%\xymatrix{
% M \ar[r] \ar[d] & M_\beta(X) \ar[d] \\
%  U \ar[r] & C}
\end{align*}
It is the disjoint union of $M_{\beta'}(X_c)$ over all $\beta'\in H_2(X_c)$
whose image in $H_2(X)$ is $\beta$, and we have
\begin{equation}\label{equ on kz2}\int_{[M_\beta(X_c)]^{\mathrm{vir}}}1=\sum_{\begin{subarray}{c}\beta' \in H_2(X_c, \mathbb{Z}) \\ i_*\beta'=\beta\end{subarray}}\int_{[M_{\beta'}(X_c)]^{\mathrm{vir}}}1,\end{equation}
where $i: X_c\to X$ is the natural closed immersion.

\subsection{Main results} 
Our main result relates \eqref{equ on kz2} to \eqref{intro equ on tau inv}. 
\begin{thm}
[Theorem \ref{222}]
\label{intro 222}
Let $c\in U$, $L=\pi^*\cO_C(c)$ and $\gamma \in H^2(X, \bQ)$. If there exists an orientation of $M_\beta(X)$ that restricts to the canonical orientation
of $M_\beta(X_c)$, then with that orientation, the following identities hold: 
$$ 
\langle\tau_0(\gamma \cup c_1(L) )\rangle_{X,\beta}=\int_\beta  \gamma\cdot \langle\tau_1( c_1(L) )\rangle_{X,\beta}  =\int_\beta  \gamma\cdot \int_{[M_\beta(X_c)]^{\mathrm{vir}}}1,
$$
\[
%\langle\tau_1(c_1(L)) \rangle_{\beta}=\deg \bfv{M_\beta(X_c)}
%, \quad
\frac{1}{12}\langle\tau_0(c_2(X)) \rangle_{X,\beta}+\langle\tau_2(1_X) \rangle_{X,\beta}=0.
\]
In particular, this confirms \cite[Conjecture 2.4]{CMT1} and \cite[Eqn.~in~\S 1.7]{CT2}. 
\end{thm}
We briefly discuss the proof of Theorem \ref{intro 222}. 
%The proof of this theorem follows from corresponding result on comparing virtual classes. 
Let $M$ be the preimage
of $U$ under $M_\beta(\pi) \to C$, i.e.,
\begin{align*}
\xymatrix{
M \ar[r]^{ } \ar[d]  \ar@{}[dr]%|{\Box}  
	& M_\beta(\pi) \ar[d] \\
U \ar[r]%^{\iota} 
	&  C}
%\xymatrix{
% M \ar[r] \ar[d] & M_\beta(X) \ar[d] \\
%  U \ar[r] & C}
\end{align*} is a Cartesian diagram. 
As an open subscheme of $M_\beta(X)$, $M$ has an induced obstruction theory which has a \textit{canonical choice of orientation} (Proposition \ref{prop on can ori}).
For $c\in U$, the pullback of this canonical orientation to $M_\beta(X_c)$ is called the \textit{canonical orientation} of $M_\beta(X_c)$. 

From another point of view, $M$ is the moduli space of stable sheaves on the $U$-family of 
Calabi-Yau 3-folds $\{X_c\}_{c\in U}$, which has a virtual class $[M]^{\mathrm{vir}}$ \cite{BF1,Tho}.
For $c\in U$, the Gysin pullback of $[M]^{\mathrm{vir}}$ along $\{c\} \to U$ is clearly $[M_\beta(X_c)]^{\mathrm{vir}}$. 
\begin{thm} 
[Corollary \ref{DT43} which is based on Theorem~\ref{OTBF}]
Take $c\in U$ and assume there exists an orientation of $M_\beta(X)$ that restricts to the canonical orientation of $M_\beta(X_c)$. Then with that orientation, we have
\[
i_c^![M_{\beta}(X)]^{\mathrm{vir}}= [M_\beta(X_c)]^{\mathrm{vir}}.
\] 
\end{thm}
Combining the above equality with a calculation of the insertion for $L=\pi^*\cO_C(c)$ (Lemma~\ref{lem on tau}), 
we obtain 
$$\int_{[M_{\beta}(X)]^{\mathrm{vir}}}c_1(L^{[1]})=\int_{[M_{\beta}(X)]^{\mathrm{vir}}}\cO_C(c)_{|_{M_\beta(X)}}=\int_{[M_\beta(X_c)]^{\mathrm{vir}}}1, $$
where $L^{[1]}$ is the tautological complex associated with $L$ \eqref{def of tau cpx}. 
Finally by a calculation using Grothendieck-Riemann-Roch, we obtain Theorem \ref{intro 222}.

We remark that Park \cite{Par} has established the functorial property of $\DT_4$ virtual class which is a powerful tool 
for many applications, e.g., proving conjectures of \cite{CK1, CKM1, CMT1, CT1} in many cases. One can also prove our theorem using his work.

\section{Preliminaries}\label{sect on pre}
\subsection{Notation, terminology, and convention}
It is convenient to assume schemes are separated, finite type over $\bC$ starting from the next subsection. In particular, schemes are separated Noetherian and maps between them are quasi-compact and quasi-separated (qcqs).

\subsubsection{Derived categories and derived functors}
For a scheme X, denote 
$\cO_X$ or $\cO$ its structure sheaf, and 
$\rmd(X)$, $\perf(X)$, $\coh(X)$, and $\qcoh(X)$ respectively, the derived category of $\cO_X$-modules,
its subcategory of perfect complexes, of complexes with coherent or quasi-coherent cohomology sheaves. As usual, superscripts including $b, -,+$ indicates
the range of nonzero cohomology sheaves, and $\tau$ denotes the truncation functor on $\rmd(X)$, with superscripts
indicating the range of truncation.

We have derived functors 
\[
\otimes, \HOM\colon \rmd(X)\times \rmd(X) \to \rmd(X),\]
with which $\rmd(X)$ is symmetric monoidal.
The adjunction between
$(-)\otimes E $ and $\HOM(E,-)$ for any object $E \in D(X)$ will be referred to
as $\tensorhom$ adjunction. 
Note that for $X$ locally Noetherian, $\HOM$ maps  $\coh^{-}(X) \times \qcoh^+(X)$ to $\qcoh^+(X)$.  
For $E\in \rmd(X)$, its dual $\HOM(E,\cO_X)$ is denoted $E^\vee$. 

For $E,F\in \rmd(X)$, we have the \textit{evaluation map} 
\[\ev\colon E \otimes \HOM(E,F) \to F.\]

 Let  $f \colon X \to Y$ be a map between qcqs schemes.
 %f qcqs suffice
We have adjunctions 
\[
f^*\colon \rmd(Y) \rightleftarrows \rmd(X)\colon\!f_*.
\] 
We will only use the right adjoint of $f_*\colon \qcoh(X) \to \qcoh(Y)$
(\cite[\href{https://stacks.math.columbia.edu/tag/0A9E}{0A9E}]{stacks})
when $f$ is proper, it will be denoted $f^!$, and $\omega_f^\bullet=f^!\cO_Y$ denotes the dualizing complex.

We will call the unit or counit map of an adjunction simply the adjunction map.
The composition of the functor $f_*$ and $\HOM$ on $\rmd(X)$ will also be written as $\HOM_f$.
When working with sheaves written in
calligraphic font, we will use
$\mathbf{R}f_*$, $\mathbf{L}f^*$, $\rhom(,)$ to emphasis functors that are derived.
The identity map of an object $E$ is denoted $\id_E$.

\subsubsection{Cotangent complex and obstruction theory}
 Let $f\colon X \to Y$ be a map of schemes, the cotangent complex of $f$
 is denoted $\bL_f \in \qcoh(X)$(\cite[\href{https://stacks.math.columbia.edu/tag/08T2}{08T2}]{stacks}),
  and $\bL^{\geqslant -1}_f$ the truncated cotangent complex. 
 If $Y$ is locally Noetherian and $f$ locally of finite type, then $\bL_f \in \coh(X)$, 
and an obstruction theory for $f$ is a map
$
\phi\colon E \to \bL_f 
$ or $E \to \bL^{\geqslant -1}_f$ in $\coh(X)$ whose cone belongs to $\coh^{\leqslant -2}(X)$, this is 
equivalent to the condition that $h^0(\phi)$ is an isomorphism, $h^{-1}(\phi)$ a surjection, and $E \in \coh^{\leqslant 0}(X)$.

\subsubsection{Nondegenerate 2-forms and maximal isotropic maps} \label{isocplx}
Let $Y$ be a scheme, and $E \in \perf(Y)$. A $n$-\textit{shifted 2-form} on $E^\vee$ is
an alternating map 
\[
\alpha \colon E^\vee \otimes E^\vee \to \cO_Y[n] \]
in $\perf(Y)$,
or equivalently, it is a map
\[
\vartheta\colon E^\vee[-n] \to E
\]
such that
$
\vartheta^\vee[-n]=-\vartheta
$ under the isomorphism 
\[
(E^\vee[-n])^\vee[-n] \simeq ((E^\vee)^\vee[n]))[-n]\simeq E.
\]
The 2-form $\vartheta$ is \textit{non-degenerate} if it is an isomorphism.
%%%
Note that the alternating map $\alpha$ can be viewed as a symmetric map
\[
E^\vee[1] \otimes E^\vee[1] \simeq (E^\vee \otimes E)^\vee[2] \xrightarrow{\alpha[2]} \cO_Y[n+2].
\]

When the 2-form $\vartheta$ is non-degenerate,
a map $\sigma\colon E \to F$ in $\perf(Y)$ is \textit{isotropic} if the composition
\[
F^\vee[-n] \xrightarrow{\vartheta\circ \sigma^\vee[-n]} E \xrightarrow{\sigma} F
\] is zero in the derived category. We will say $\sigma$ is \textit{maximal isotropic} if the induced map from $F^\vee[-n]$
to the cocone of $\sigma$ is an isomorphism, and will
refer to $\sigma$ by saying  $F$ is a maximal isotropic complex of $E$.

\begin{rmk}

When $Y=\spec\,\bC$, $n=0$, and $E$ a finite dimensional complex vector space, a non-degenerate 2-form is
a symplectic structure on $E$.
%A non-degenerate $n$-shifted 2-form on $E$ is the same as
%a non-degenerate $(-n)$-shifted 2-form on $E^\vee$.
The terminology ``maximal isotropic"  introduced above is not standard. 
\end{rmk}

\subsubsection{}

The map $\pr_X$ denotes the projection map $X \times_S Y \to X$ for schemes $X,Y$ over $S$.

\subsection{A self-dual distinguished triangle}

Let $i\colon Z \to X$ be the inclusion of an effective Cartier divisor
%(\cite[\href{https://stacks.math.columbia.edu/tag/01WR}{01WR}]{stacks})
into a scheme $X$.
%No finiteness assumptions is needed on $X$.
The right adjoint $i^!$ of 
$$i_*\colon \rmd(Z) \to \rmd(X)$$
is given by (\cite[\href{https://stacks.math.columbia.edu/tag/0A76}{0A76}]{stacks}):
$$\rhom(\cO_Z, -),$$ 
which satisfies (\cite[\href{https://stacks.math.columbia.edu/tag/0AA4}{0AA4}]{stacks}):  
$$\rhom(\cO_Z, -) \simeq \mathbf{L}i^*(-)\otimes \cN[-1], $$
where $\cN$ is the normal bundle of $i$. Note that $\cN[-1]$ is the dualizing complex $\omega_i^\bullet$ defined as $i^!\cO_X$.

For any $F \in \rmd(Z)$, the adjunction map 
$$F \to i^!i_*F\simeq  i^*i_*F \otimes \cN[-1]$$ determines
by tensoring with $\cN^\vee[1]$ the map
$F\otimes \cN^\vee[1]\to i^*i_*F $ . We also have the adjunction map 
$i^*i_*F \to F$.  It is argued using Fourier-Mukai kernels (\cite[Lemma 3.3]{BO},~\cite[Corollary. 11.4 (ii)]{Huy}) that 
%when $X$ is (smooth) separated and Noetherian
\begin{equation} \label{dist}
F\otimes \cN^\vee[1]\to i^*i_*F \to F
\end{equation} is a distinguished triangle for $F \in \coh^b(Z)$.

\begin{rmk} 
%Assume $X$ is a separated Noetherian scheme.  
By \cite[Lemma 2.2]{HT}, on $Z \times Z$ there exists an exact sequence 
\begin{equation}\label{FM}
0 \to {\Delta}_*\cN^\vee \to {\cI_X}_{|_{Z\times Z}} \to \cO_{Z\times Z} 
\to {\Delta}_*\cO_Z \to 0,
\end{equation} where $\Delta \colon Z \to Z \times Z$ is the diagonal emedding, and  $\cI_{X}$ is the ideal sheaf of the diagonal $X$ in $X \times X$.
The complex consisting of the middle two terms ${\cI_X}_{|_{Z\times Z}} \to \cO_{Z\times Z}$  with $\cO_{Z\times Z}$
at degree $0$ is the Fourier-Mukai kernel of 
$i^*i_*$, and the exact sequence viewed as a triangle of Fourier-Mukai kernels induces for any $F\in \qcoh(X)$
a distinguished triangle
\[
F\otimes \cN^\vee[1]\to i^*i_*F \to F,
\] which is isomorphic to \eqref{dist}.
% up to automorphisms of the kernels $\Delta_*N^\vee[1]$ and $\Delta_*\cO_Z$.
For $F$ a finite type quasi-coherent sheaf, a transparent proof that \eqref{dist} is distinguished is given in Lemma \ref{triangle}. See Remark \ref{qcs} for the general case.
\end{rmk}

%\begin{rmk}
%It follows from the previous remark and arguments in \cite[Section 2.4]{HT}  that, if $X$ embeds into a smooth and separated Noetherian scheme, then the map $F \to F\otimes \cN^\vee[2]$ in the triangle \eqref{dist}  is up to isomorphisms induced by the truncated Atiyah class $F \to F\otimes \bL_Z^{\geqslant -1}[1]$ and the map between cotangent complexes $\bL_Z \to \bL_{i}\simeq \cN^\vee[1]$. When
%$X$ is fibered over a curve and $Z \to X$ the inclusion of a fiber, the map $F \to F\otimes \cN^\vee[2]$
%is the obstruction to first order deformation of $F$ off the fiber $Z$(See, e.g., \cite[Theorem 3.45]{Tho}).
%\end{rmk}

Applying $i_*\HOM(-, F)$ to \eqref{dist},
and using the isomorphism (\ref{homlst}): 
$$\HOM(i_*F,i_*F)\simeq i_*\HOM(i^*i_*F, F), $$ 
we obtain a distinguished triangle
\[
i_*\HOM(F,F) \to \HOM(i_*F,i_*F)\to i_*\HOM(F, F\otimes \cN[-1]),
\]
for the last term we used the canonical isomorphism
$$\HOM(F\otimes \cN^\vee[1], F)  \simeq \HOM(F, F\otimes \cN[-1]). $$

\begin{prop} \label{selfdual}
Let $i\colon Z \to X$ be the inclusion of an effective Cartier divisor into a separated Noetherian scheme $X$.
If $F \in \perf(Z)$, then the triangle 
\[
i_*\HOM(F,F) \to \HOM(i_*F,i_*F)\to i_*\HOM(F, F\otimes \cN[-1])
\]
is self-dual.
\end{prop}

\begin{proof}
We will identify the dual of $\HOM(i_*F,i_*F)\to i_*\HOM(F, F\otimes \cN[-1])$ with 
$i_*\HOM(F,F) \to \HOM(i_*F,i_*F)$.
Note that $i_*$ preserves perfect objects by
(\cite[\href{https://stacks.math.columbia.edu/tag/068C}{068C},
\href{https://stacks.math.columbia.edu/tag/0B6G}{0B6G}]{stacks}).

Let $K,L \in \perf(Z)$.
Applying 
$i_*\HOM(-,L)$ to $i^*i_*K \to K$, we get 
\[
i_*\HOM(K,L) \to i_*\HOM(i^*i_*K, L)\simeq \HOM(i_*K, i_*L).
\] This is the map $\hom_i$\eqref{homf} by Lemma \ref{dual1}.

Applying $(i_*\HOM(-,K))^\vee$ to $L\otimes \cN^\vee[1]\to i^*i_*L$,
using the following isomorphism of functors  on $\perf(Z)$
\[
(i_*\HOM(-,K))^\vee \simeq i_*\HOM(\HOM(-,K), \omega_i^\bullet) \simeq i_*\HOM(K, (-)\otimes \omega_i^\bullet),
\]
we get 
\[
i_*\HOM(K, L) \to i_*\HOM(K, i^!i_*L),
\]which is also induced by applying $i_*\HOM(K,-)$ to $L \to i^!i_*L$, and is identified with
$i_*\HOM(K,L) \to \HOM(i_*K, i_*L)$ by Lemma \ref{dual2}.
The proposition follows by setting $K=L=F$.

The identification is given by some commutative diagram
\[
\xymatrix{
(i_*\HOM(F, F\otimes \cN[-1]))^\vee \ar[r]\ar[d]
	&   \HOM(i_*F,i_*F)^\vee\ar[d] \\
 i_*\HOM(F,F) \ar[r] 
 	& \HOM(i_*F,i_*F).
}
\]
Note that the vertical isomorphism
$\HOM(i_*F,i_*F)^\vee \to \HOM(i_*F,i_*F)$ agrees with the obvious one
 by setting $K=F, G=i_*F, H=\cO$ in Lemma \ref{dual3}.
\end{proof}

 \subsection{Relative moduli spaces of sheaves} \label{moduli}

Let 
$$\pi \colon X \to S$$ be a projective map 
between complex projective schemes with connected fibers, and 
$\beta\in H_2(X,\bZ)$ a curve class satisfying $\pi_*\beta=0$.
Denote $M_\beta(\pi)$ the moduli scheme of 1-dimensional 
semi-stable sheaves $\cF$ on the fibers of $\pi$ with $[\cF]=\beta$ and $\chi(\cF)=1$ (e.g.,~\cite[Theorem 4.3.7]{HL}). 
Since $\chi(\cF)=1$, semi-stability coincides with stability and $M_\beta(\pi)$ is a projective fine moduli scheme over $S$ with a universal family. 

\begin{rmk}
The condition $[\cF]=\beta$ is interpreted as follows.
For any scheme $T$ over $S$, and a family of sheaves $\cE$ over $T \times_S X$, by pushing forward along the closed immersion $T \times_S X \to T \times X$, we view $\cE$ as a family over $T \times X$, and require 
$[\cE_{|_{t \times X}}]=\beta$ for any closed point $t\in T$.
For a coherent sheaf $\cF$ on $X$, $[\cF]$ denotes the cycle class of the 1-cycle $Z_1(\cF)$
(see \cite[Example~18.3.11]{Ful}).
\end{rmk}

%\begin{rmk}
%Note that stability is independent of the choice of polarization, and being semi-stable implies being stable. (See e.g.,\cite{Katz}.) In fact, for any 
%$\pi$-ample line bundle $\cL$, the Hilbert Polynomial of $\cF$ is $\chi(\cF\otimes \cL^{\otimes n})=(\int_\beta c_1(\cL)) n+\chi(\cF)$,
%and it follows that a pure 1-dimensional sheaf $\cF$ is semi-stable if and only if for any proper subsheaf $\cF'$, $\chi(\cF')\leqslant 0$. 
%\end{rmk}

\begin{rmk} \label{norm}
The existence of universal family on $M_\beta(\pi)$ follows from arguments in \cite[Section~4.6]{HL}.
Let $\sM_\beta(\pi)$ be the corresponding moduli stack, it is a $\mathbb{C}^*$-gerbe over $M_\beta(\pi)$. 
Quasi-coherent sheaves on a $\bC^*$-gerbe decompose into direct summands indexed by integer weights(see e.g., \cite[Prop. 2.2.1.6]{Lie}).
Denote $\sE$ the universal sheaf  on $\sM_\beta(\pi)\times_S X$ which has weight 1, and 
$\pr \colon\sM_\beta(\pi) \times_S X \to  \sM_\beta(\pi)$
the projection map.
As the line bundle $\det(\mathbf{R}\pr_*\sE)$ has weight $\chi(\sE)$,
$\sE \otimes \pr^*(\det(\mathbf{R}\pr_*\sE))^{\vee}$ has weight 0 and descends to a universal sheaf on $M_\beta(\pi) \times_S X$, it is the  “normalized” one in \cite[(1.7)]{CT2}. 
\end{rmk} 
% $\pr_*\sE$ is perfect

Similarly, denote $M_\beta(X)$ the projective fine moduli scheme of 1-dimensional 
semi-stable sheaves $\cF$ on $X$ with $[\cF]=\beta$ and $\chi(\cF)=1$. We will identify the moduli space $M_\beta(X)$ with $M_\beta(\pi)$.

Consider the normalized universal sheaf
on $M_\beta(\pi)\times_S X$, pushing forward along the closed immersion 
$M_\beta(\pi)\times_S X \to M_\beta(\pi)\times X$ produces a family over $M_\beta(\pi)\times X$, and induces a map $\epsilon\colon M_\beta(\pi) \to M_\beta(X)$.

\begin{prop} \label{relabs}
The map $\epsilon \colon M_\beta(\pi) \to M_\beta(X)$
is an isomorphism. 
\end{prop}

\begin{proof}
In the following, $\HOM(\cG,\cG)$ and $p_*$ are not derived.

Let $T$ be a finite type scheme, and $\cG$ a family of 1-dimensional stable sheaves on $X \times T$ that determines a map
$T \to M_\beta(X)$.
We will use the rigidity lemma below to show the map 
\[
\Supp(\cG) \hookrightarrow X \times T  \xrightarrow{\pr_X} X \xrightarrow{\pi} S
\]
factors through the map $p\colon \Supp(\cG) \hookrightarrow X \times T  \xrightarrow{\pr_T} T$.
Here $\Supp(\cG)$ denotes the support of $\cG$, it is a closed subscheme of $X \times T$ with ideal sheaf
the kernel of $
\cO_{X\times T} \xrightarrow{}\HOM(\cG,\cG)
$. Clearly $\Supp(\cG) \to T$ is projective, hence closed.

The natural map $\cO_{X\times T} \xrightarrow{}\HOM(\cG,\cG)$ factors as 
\[
\cO_{X\times T} \xrightarrow{} \cO_{\Supp(\cG)}\hookrightarrow \HOM(\cG,\cG).
\] Note that $ \cO_{\Supp(\cG)}\hookrightarrow \HOM(\cG,\cG)$ is injective.
Pushing forward along $\pr_T$, we obtain
\[
{\pr_T}_*\cO_{X\times T} \to p_*\cO_{\Supp(\cG)} \hookrightarrow {\pr_T}_*\HOM(\cG,\cG).
\]
As the composition $\cO_T \to {\pr_T}_*\cO_{X\times T} \xrightarrow{} {\pr_T}_*\HOM(\cG,\cG)$ is an isomorphism, we conclude
\[
\cO_T \simeq p_*\cO_{\Supp(\cG)}.
\]
For any closed point $t \in T$, its preimage in $\Supp(\cG)$ as a set is the 
underlying set of $\Supp(\cG_{|_{X \times t}})$ (\cite[\href{https://stacks.math.columbia.edu/tag/056J}{056J}]{stacks}).
As $S$ has ample line bundles and $\pi_*\beta=0$, any 
irreducible component of $\Supp(\cG_{|_{X \times t}})$, which has dimension 1,  is mapped to a closed point of $S$.
As $\cG_{|_{X \times t}}$ is simple, $\Supp(\cG_{|_{X \times t}})$ can not be mapped to more than one point of $S$.
Now we can apply the lemma below and finish the proof.
\end{proof}

The version of rigidity lemma we use is the following, 
and we include a proof for convenience of the reader.
\begin{lemma}
Let $f\colon Y \to T$ and $g\colon Y \to S$ be maps between arbitrary schemes. 
Assume $f$ is a closed map, $\cO_T\simeq f_*\cO_Y$, and $t_0 \in T$ is a point such that  
$g$ maps $f^{-1}(t_0)$ to a single point $s$ of $S$. Then there exists an open neighborhood
$U$ of $t_0$ in $T$, such that the restriction of $g$ to $f^{-1}U$ factors uniquely through $f_{|_{f^{-1}U}}$.
It follows that $g$ factors through $f$ if $T$ is locally Noetherian and for any closed point $t \in T$,
$g(f^{-1}(t))$ is a single point.
\end{lemma}

\begin{proof}
Let $V$ be an open affine neighborhood of $s$ in $S$, then $g^{-1}(S-V)$ is closed and disjoint from 
$f^{-1}(t_0)$.
Let $U$ be the complement of $f(g^{-1}(S-V))$, then $U$ is an open neighborhood of $t_0$, and $f^{-1}U $ is contained in $g^{-1}V$, or $g$ maps $f^{-1}U$ into $V$.

Denote $f_U$ and $g_U$ the restriction of $f$ and $g$ to $f^{-1}U$ respectively.
As $V$ is affine, the map $g_U\colon f^{-1}U \to V$ is determined by 
\begin{equation} \label{map}
\Gamma(V, \cO_V) \to \Gamma(f^{-1}U, \cO_Y),
\end{equation}
and $g_U$ factors through $f_U$ iff
\eqref{map} factors through
\[
\Gamma(U, \cO_T)\xrightarrow{} (f^{-1}U, \cO_Y).
\] This certainly holds as the map is an isomorphism by $\cO_T \simeq f_*\cO_Y$.

Consider the subset of points $t \in T$ such that $g(f^{-1}(t))$ is a single point, then
it is open. 
When $T$ is locally Noetherian, its complement must be empty, otherwise it contains a closed point (\cite[\href{https://stacks.math.columbia.edu/tag/02IL}{02IL}]{stacks}).
\end{proof}

 \section{Compatibility between obstruction theories}\label{sect on pre2}
 
 In this section, we recall an alternative formulation of the
standard obstruction theory used in sheaf-counting theories (Lemma \ref{alt}),
and use it to verify a compatibility result between obstruction theories (Proposition \ref{phipf}).
 
 \subsection{Equivalent formulations of obstruction theory maps}
 \label{phi}
 
Let $\pi\colon X \to S$ be a proper and flat map between Noetherian schemes, 
$T$ a scheme over $S$, and $f\colon Y \to T$ the base change of $\pi$ to $T$. As $f$ is proper and flat, 
 $f_*$ maps perfect objects to perfect objects
 (\cite[\href{https://stacks.math.columbia.edu/tag/0B6F}{0B6F}]{stacks}),
 and the right adjoint $f^!$ of $f_*$
satisfies
$$f^!(-)\simeq f^* (-) \otimes \omega_f^\bullet, $$
where $\omega_f^\bullet =f^!\cO_T$ is the \textit{dualizing complex} of $f$.
(\cite[\href{https://stacks.math.columbia.edu/tag/0E4K}{0E4K}]{stacks}).

Let $E \in \perf(Y)$, then $\HOM(E,E)$ and $f_*\HOM(E,E)$ are perfect. Consider the \textit{relative Atiyah class} of $E$
\[\at(E)\colon E \to E\otimes \bL_{Y/X}[1],
\] 
where $ \bL_{Y/X} \in \qcoh(Y)$ denotes
cotangent complex of the projection map $\pr_X\colon Y \to X$. 
Using the map $\at(E)$,  we have the map 
\[
\psi_E \colon f_*(\HOM(E,E)\otimes \omega_f^\bullet)[-1] \to \bL_{T/S}
\] 
in the construction of obstruction theories in sheaf-counting theories.
Recall $\at(E)$ determines by adjunction and shift
\[
\HOM(E,E)[-1] \to \bL_{Y/X} \simeq f^*\bL_{T/S},
\]
tensoring with $\omega_f^\bullet$, we obtain
$\HOM(E,E)[-1]\otimes \omega_f^\bullet \to f^!(\bL_{T/S})$, then
by $\usk{f}$ adjunction, we arrive at $\psi_E$:
\begin{equation}
  \label{psie} f_*(\HOM(E,E)\otimes \omega_f^\bullet)[-1] \to \bL_{T/S}.
\end{equation}

An alternative formulation is given by as follows.
As $E$ is perfect, from $\at(E)$ we obtain
the map 
\begin{equation} \label{ate} 
\At(E)\colon \cO_Y \to \HOM(E,E)[1]\otimes \bL_{Y/X}.
\end{equation}
As $ f^*\cO_T \simeq \cO_Y$, by $\lst{f}$ adjunction, 
\[
\cO_Y \to \HOM(E,E)[1]\otimes \bL_{Y/X}
\] corresponds to 
\[
\cO_T \to f_*(\HOM(E,E)[1]\otimes f^*\bL_{T/S})\simeq
 f_*\HOM(E,E)[1]\otimes \bL_{T/S}
\] 
where the isomorphism is given by the projection formula (\ref{projfor}).
As $f_*\HOM(E,E)[1]$ is perfect, we obtain
\begin{equation}
 \label{phie}
(f_*\HOM(E,E)[1])^\vee \to \bL_{T/S},
\end{equation}
henceforth denoted $\phi_E$.
Abusing notation, we will also use $\phi_E$ later to denote maps from 
some truncation of $(f_*\HOM(E,E)[1])^\vee$ to  the truncated cotangent complex
$\bL_{T/S}^{\geqslant -1}$.

\begin{rmk}
The map $\phi_E$ is the same as the map $\mathrm{At}_E$ in \cite[Section 1.4]{Ku}.
\end{rmk}

As a special case of Lemma \ref{pots}, we have
\begin{lemma} \label{alt}
Let $f\colon Y \to T$ be a proper flat map and $E \in \perf(Y)$.
Under the duality isomorphism
\[
(f_*\HOM(E,E)[1])^\vee \simeq  f_*(\HOM(E,E)\otimes \omega_f^\bullet)[-1],
\]
the maps $\phi_E$ \eqref{phie} and $\psi_E$ \eqref{psie} are the same.
\end{lemma}

%\begin{rmk}
%To check compatibility between obstruction theories, $\phi_E$ seems to be easier to work with.
%\end{rmk}

\subsection{Functoriality of $\phi$ under pushforwards}  \label{functoriality}

As before $\pi\colon X \to S$ denotes a proper flat map.
Let $M$ be a scheme over $S$, and $F \in \perf(M\times_SX)$. 
Denote $p\colon M\times_S X \to M$ and $q\colon M\times X \to M$
the projection maps to $M$. Denote $j \colon M\times_S X \to M \times X$ the closed immersion, and assume $G=j_*F$ is perfect.

Using the Atiyah classes 
\[\at(F)\colon F \to F\otimes \bL_{M\times_SX/X}[1],\quad \at(G)\colon G \to G\otimes \bL_{M\times X/X}[1], \]
we obtain 
\[\phi_F\colon (\HOM_p(F,F)[1])^\vee \to \bL_{M/S}, \quad \phi_G\colon (\HOM_q(G,G)[1])^\vee \to \bL_M,\] 
and we will show the compatibility between $\phi_G$ and $\phi_F$.

\begin{lemma}
If $S$ is smooth and separated, then $G=j_*F$ is perfect.
\end{lemma}

\begin{proof}

Consider the Cartesian diagram
\[
\xymatrix{
M \times_SX \ar[r]^j\ar[d] 
	& M\times X\ar[d]\\
S \ar[r]^{\Delta_S}
	& S \times S.
}
\]
As $S$ is smooth and separated, the diagonal map $\Delta_S\colon S \to S\times S$ is a regular closed immersion.
It follows from flatness of $\pi$ that the diagram is tor-independent, as $M\times_S X$ is the derived fiber product. From tor-independence of the Cartesian diagram, we conclude $j$ is a regular closed immersion.
In particular, $j$ is perfect proper so that $j_*$ preserves perfect objects.
\end{proof}

\begin{prop}Assume $G$ is perfect. \label{phipf}
The diagram
\[
\xymatrix{
(\HOM_q(G,G)[1])^\vee \ar[r]\ar[d]^{\phi_G} 
	&  (\HOM_p(F,F)[1])^\vee \ar[d]^{\phi_F}\\
\bL_{M} \ar[r] & \bL_{M/S}.
}
\] is commutative. The vertical arrows are as indicated $\phi_G, \phi_F$.
The upper horizontal arrow is induced by
$\hom_j\colon  j_*\HOM(F,F) \to \HOM(j_*F, j_*F)$(\ref{homf}), and
 lower horizontal arrow the natural map between cotangent complexes.
 Note that for $\phi_G$, the map $q$ comes from the trivial family $X \to \spec\bC$.
\end{prop}

\begin{proof}
From the adjunction map $j^*G\simeq j^*j_*F \to F$ and by functoriality and compatibility with respect to pullbacks of Atiyah classes (see e.g., \cite[1.1]{Ku}), we 
obtain the commutative diagram
\[
\xymatrix{
j^*G \ar[r]^-{j^*\at(G)} \ar[d]
 	& j^*G\otimes j^*\bL_{M\times X/X}[1] \simeq j^*(G\otimes q^*\bL_M)[1]\ar[d]\\
F\ar[r]^-{\at(F)}
 	& F \otimes \bL_{M\times_SX/X}[1]\simeq F \otimes p^*\bL_{M/S}[1].
}
\]
By $\lst{j}$ adjunction, this corresponds to the commutative diagram
\begin{equation} \label{comm}
\xymatrix{
G \ar[r]^-{\at(G)}\ar[d]
 & G\otimes q^*\bL_{M}[1]\ar[d]\\
j_*F\ar[r]^-{j_*\at(F)}
 & j_*F \otimes q^*\bL_{M/S}[1],
}
\end{equation}
where we used the projection formula 
\[
j_*F \otimes q^*\bL_{M/S} \simeq j_*(F \otimes j^*q^*\bL_{M/S}) \simeq j_*(F \otimes p^*\bL_{M/S}).
\] 
%and we see the composition of $\at(G)$ and 
%$G\otimes q^*\bL_M[1]\to G\otimes q^*\bL_{M/S}[1]$ induced by
%$\bL_M \to \bL_{M/S}$ is identified with $j_*\at(F)$.

Consider the commutative diagram
\[
\xymatrix{
j_*F\otimes \cO_{M\times X}\ar[d] \\
j_*F\otimes j_*\cO_{M\times_S X} \ar[r]^-{\id\otimes j_*\At(F)} \ar[d]^{\otimes_j\eqref{otimesf}}
	& j_*F \otimes j_*(\HOM(F,F)[1]\otimes p^*\bL_{M/S}) \ar[d]^{\otimes_j}
		&\\
j_*(F\otimes \cO_{M\times_S X}) \ar[r]^-{j_*\circ\id\otimes \At(F)}
	& j_*(F \otimes \HOM(F,F)[1]\otimes p^*\bL_{M/S}) \ar[r]^-{j_*\circ \ev\otimes \id}
		&j_*(F \otimes  p^*\bL_{M/S}[1]),
}
\] where the unspecified map is defined by tensoring $j_*F$ with $\cO_{M\times X} \to j_*\cO_{M\times_SX}$.
Note that the composition of the two arrows in the botton line is $j_*\at(F)$.
As the left vertical composition
\[
j_*F\otimes \cO _{M\times X}\to j_*(F\otimes \cO_{M\times_S X})
\] is the identity map of $j_*F$ under $j_*F\simeq j_*F\otimes \cO$ and $F \simeq F\otimes \cO$,
%\footnote{By $\lst{j}$ adjunction, the identification with  $\id_{j_*F}$ comes down to compatibility between the units $\cO_{M\times X}$ and $\cO_{M\times_SX}$ with respect
%to the symmetric monoidal functor $j^*$.}
the zigzag 
from $j_*F\otimes \cO_{M\times X}$ to $j_*(F\otimes p^*\bL_{M/S}[1])$ is
identified with $j_*\at(F)$. It follows from the adjunction between $j_*F\otimes (-)$ and $\HOM(j_*F, -)$ on $\rmd(M\times X)$, the compatibility between $j_*\at(F)$ and $\at(G)$ \eqref{comm}, 
and the compatibility between $\otimes_j$ and the projection formula that  the diagram
\[
%\begin{equation} 
\xymatrix{
\cO_{M\times X} \ar[r] \ar[d]^-{\At(G)\eqref{ate}} 
	& j_*\HOM(F,F)[1]\otimes q^*\bL_{M/S} \ar[d]^-{\hom_j [1]\otimes \id}  \\
 \HOM(G,G)[1] \otimes q^*\bL_M\ar[r]
	&  \HOM(G,G)[1] \otimes q^*\bL_{M/S}.
}
\]
%\end{equation}  
is commutative, where
the horizontal arrow 
\[
\cO_{M\times X} \to j_*\HOM(F,F)[1]\otimes q^*\bL_{M/S}\simeq j_*(\HOM(
F,F)[1] \otimes p^*\bL_{M/S})
\]
corresponds under $\lst{j}$ adjunction to $\At(F)$\eqref{ate} . 
Then by $\lst{q}$ adjunction, we see
 \[
\xymatrix{
\cO_M \ar[r]\ar[d] 
	& q_*( j_*\HOM(F,F)[1]\otimes q^*\bL_{M/S}) \ar[d] 
		&  p_*\HOM(F,F)[1]\otimes \bL_{M/S} \ar[l]^-{\simeq}_-{\mathrm{pf}} \\
q_*(\HOM(G,G)[1] \otimes q^*\bL_M)	\ar[r]	
	& q_*(\HOM(G,G)[1] \otimes q^*\bL_{M/S})
		& \\
\HOM_q(G,G)[1]\otimes \bL_M \ar[u]^-{\simeq}_-{\mathrm{pf}} \ar[r]
	& \ar[u]^-{\simeq}_-{\mathrm{pf}} \HOM_q(G,G)[1]\otimes \bL_{M/S},
		&
}
\] is commutative. Now the proposition follows from the construction of $\phi_G$, $\phi_F$,
and the fact that $\lst{p}$ adjunction is the composition of $\lst{j}$ and $\lst{q}$ adjunctions.
\end{proof}

\subsection{Functoriality of $\phi$ under pullbacks}

Let $X$ be a scheme considered as a trivial family $X \to \spec\bC$, and $\iota\colon M \to N$  a map between schemes over $\bC$.
Denote $q\colon M\times X \to M $ and $r\colon N\times X \to N$ the pullback families. Consider the Cartesian diagram
\[
\xymatrix{
M\times X \ar[r]^-{\kappa} \ar[d]^-{q}
	& N \times X\ar[d]^-{r} \\
M \ar[r]^-{\iota} 
	& N
} 
\] where $\kappa=(\iota, \id_X)$.
\begin{prop} \label{phipb}
 Let  $E \in \perf(N \times X)$, and $G=\kappa^*E$.
Under the isomorphism
$$\iota^*(r_*\HOM(E,E)[1])^\vee \simeq q_*(\HOM(\kappa^*E,\kappa^*E)[1])^\vee$$ induced by
$\iota^*\circ r_*\simeq q_*\circ \kappa^*\colon \qcoh(N \times X) \to \qcoh(M)$
(\cite[\href{https://stacks.math.columbia.edu/tag/08IB}{08IB}]{stacks}),
we have a commutative diagram
\[
\xymatrix{
\iota^* (\HOM_r(E,E)[1])^\vee \ar[d]^-{\iota^*\phi_E} 	\ar[r]^-{\sim} 
	&  (\HOM_q(G,G)[1])^\vee \ar[d]^-{\phi_G}\\
\iota^*\bL_N \ar[r] 
	&\bL_M,
}
\] where  bottom horizontal arrow is the map between cotangent complexes induced by $\iota$.
\end{prop}

This is straightforward to prove, and we omit the details as the result is not used in the paper.
 
 \section{DT4 virtual classes}

Let $X$ be a connected smooth complex projective CY 4-fold, and 
$M$ a quasi-projective moduli scheme of 
 Gieseker stable sheaves on $X$. Assume there is a universal
(twisted) sheaf $E$ on $M\times X$, then the Atiyah class of $E$ induces an obstruction  theory 
\begin{equation}\label{obs}
\phi_E\colon \tau^{[-2,0]}(\HOM_{\pr_M}(E,E)[1])^\vee\to \bL_M
\end{equation} for $M$ (e.g., \cite[Theorem 4.1]{HT}).
The complex $\tau^{[-2,0]}(\HOM_{\pr_M}(E,E)[1])^\vee$ is perfect
of tor-amplitude in $[-2,0]$, and it has
a non-degenerate 2-shifted 2-form induced by duality.
An orientation of $M$ is an orientation of the complex $\tau^{[-2,0]}(\HOM_{\pr_M}(E,E)[1])^\vee$,
and orientation on $M$ always exists by \cite{CGJ} (see also \cite{CL2}).
Given an orientation of $M$, Oh and Thomas \cite{OT} have produced an algebraic virtual class $\otv{M}$ 
 from the obstruction theory \eqref{obs}. 
\begin{rmk}
Let $\sM$ be a moduli stack of Gieseker stable sheaves on $X$, it is a $\bC^*$-gerbe 
over a quasi-projective moduli scheme $M$. 
Let $\sE$ be the universal sheaf on $\sM \times X$, it can be viewed as a twisted universal sheaf 
$E$ on $M$. The obstruction theory $\phi_E$ for $M$
is induced by $\phi_\sE$ for $\cM$, and the virtual classes on $\cM$ and $M$ are related by flat pullback along  $\sM \to M$ (\cite{HT,Q}).
\end{rmk}
In this section, we 
show maximal isotropic complexes of tor-amplitude in
$[-1,0]$ determine orientations
(Lemma \ref{lemmaori}) and 
allows us to identify Oh-Thomas virtual classes with Behrend-Fantechi
virtual classes (Proposition \ref{OTviaBF}). Then in Section \ref{maxiso}, under Assumption \ref{lego} we show how to construct such maximal isotropic complexes (Theorem \ref{K}).

\subsection{Orientations} 

As we will work with virtual classes defined in \cite{OT}, we follow implicitly their sign convention for orientation.

Let $Q$ be an orthogonal bundle on a scheme $M$, i.e., $Q$ a finite rank locally free sheaf
with a nondegenerate symmetric 2-form $\theta\colon Q \to Q^\vee$.
If $\Lambda \subset Q$ is a maximal isotropic subbundle, then
$\Lambda$ determines an orientation of $Q$ (e.g., \cite[Section 2]{OT}).
%\begin{rmk}
%An orientation of $Q$ is an isomorphism $o\colon \cO \to \det Q$ such that
%the composition
%\[
%\cO \xrightarrow{o} \det Q \xrightarrow{\det\theta} \det Q^\vee  \simeq (\det Q)^{\vee} \xrightarrow{o^\vee} \cO
%\] 
%is given by $(-1)^{r(r-1)/2}$, where $r$ is the rank of $Q$.
%\end{rmk}
Let $V \in \perf(M)$ be perfect of tor-amplitude in $[-2,0]$ with a nondegenerate 2-form
$\vartheta:V^\vee[2]\simeq V$.
If $V$ is represented as a self-dual complex (\cite[Equation (48)]{OT})
\[
P \xrightarrow{} Q  \xrightarrow{} P^\vee,\footnote{The dual of the bundle $P$ is denoted $P^*$ in \cite{OT}.}
\]
then orientations of $V$ correspond bijectively to orientations of $Q$
(\cite[Proposition 4.2]{OT}).

%\begin{rmk}
%An orientation of $V$ is an isomorphism $o\colon \cO \to \det V$ such that
%the composition
%\[
%\cO \xrightarrow{o} \det V \xrightarrow{{\det{\vartheta}}^{-1}} \det V^\vee[2]  \simeq (\det V^{-1}) \xrightarrow{o^\vee} \cO
%\] 
%is given by $(-1)^{r(r-1)/2}$, where $r$ the rank of $V$. 
%The complex $V$ is orientable iff $\det V \simeq \cO$.
 %if $o,o'$ are orientations, then $o'^{-1}\circ o\colon \cO \to \cO$ as an element
%in $\Gamma(M, \cO)$ corresponds to a morphism $M \to \bA^1$, which is locally constant and valued in 
%$\{-1, 1\}$.
%\end{rmk}

Assume $\sigma\colon V \to W$ is maximal isotropic,
then the distinguished triangle 
\[
W^\vee[2] \to V \to W
\]
determines a canonical isomorphism 
\[
\det V \simeq \det W \otimes \det W^\vee[2] \simeq 
\det W \otimes \det W^\vee  \simeq\cO_X.
\]
At least when $W$ is perfect of tor-amplitude in $[-1,0]$, it is easy to see this determines an orientation of $V$ (see also \cite[Remark 4.2]{CZ} for a similar discussion).
\begin{lemma} \label{lemmaori}
Let $M$ be a quasi-projective scheme and $V \in \perf(M)$ be perfect of tor-amplitude in $[-2,0]$ with a nondegenerate 2-form.
If $\sigma\colon V \to K$ is maximal isotropic and $K$ perfect of tor-amplitude in $[-1,0]$, then  
$\sigma$ determines an orientation of $V$.
\end{lemma}

\begin{proof}
By the proof of \cite[Proposition 4.1]{OT}, if we fix a two term complex of finite locally free sheaves $K^{-1} \to K^0$ presenting $K$, then we can represent $V$ as a self-dual complex
$P \to Q \to P^\vee$ and 
$\sigma$ a map between these two complexes:
\[
\xymatrix{
P \ar[r] & Q \ar[r]\ar[d] & P^\vee\ar[d] \\
           & K^{-1} \ar[r]  & K^0.
}
\]
We will show this presentation of $\sigma$ determines a maximal isotropic subbundle of $Q$.

From the distinguished triangle
$K^\vee[2] \to V \to K$, we see the cone of $\sigma\colon V \to K$ 
lies in $\qcoh^{\leqslant -2}(M)$.
As the cone is represented by the complex $P \to Q \to K^{-1}\oplus P^\vee \to K^0$, it is exact in degree $-1$ and $0$.
Let $R$ be the kernel of $K^{-1}\oplus P^\vee \to K^0$, then we have a surjective map
$Q \to R$ and an exact sequence
\[
0\to R \to K^{-1}\oplus P^\vee \to K^0 \to 0,
\] 
or a quasi-isomorphism
\[
\xymatrix{
 R \ar[r]\ar[d] & P^\vee\ar[d]          \\
K^{-1} \ar[r]  & K^0.
}
\]
Thus we have a new representative $R \to P^\vee$ of $K$, and $\sigma$ is presented as 
\[
\xymatrix{
P \ar[r]
	& Q \ar[r]\ar[d] 
		& P^\vee 	\ar[d]^{\id} \\
     &R \ar[r]     
     	& P^\vee \ar@{}[ul]|{\Box}.      
}
\]
The map $K^\vee[2] \to V$ is represented by the map
between complexes
\[
\xymatrix{
P \ar[r]  \ar[d]^{\id}  
	&		R^\vee \ar[d] 
		& \\
P \ar[r]                   
	& Q \ar[r]  \ar@{}[ul]|{\Diamond} 
		& P^\vee,       
}
\] where the square $\Diamond$ is obtained by dualizing the square $\Box$.

From the long exact sequence induced by $K^\vee[2] \to V \to K$, we see
\[
R^\vee \to Q \to R 
\] is a short exact sequence. Thus $R^\vee \subset Q$ is maximal isotropic.
\end{proof}

\subsection{Compared with Behrend-Fantechi virtual class}
\begin{prop} \label{OTviaBF}
Let $X$ be a connected smooth complex projective CY 4-fold, and 
$M$ a quasi-projective moduli scheme of Gieseker stable sheaves on $X$.
Assume the obstruction theory \eqref{obs} 
factors through a maximal isotropic complex
$K$ perfect of tor-amplitude in $[-1,0]$, 
then the induced map $K \to \bL_M$ is a
perfect obstruction theory for $M$ \cite{BF1,LT}, and 
we have $$\bfv{M}=\otv{M}.$$
Here $\bfv{M}$ denotes the virtual class determined by $K \to \bL_M$,
and $\otv{M}$ denotes the virtual class determined by \eqref{obs} and the orientation determined by $K$.
\end{prop}

\begin{proof}

As $K$ is of tor-amplitude in $[-1,0]$, $K^\vee[2] \in \qcoh^{[-2,-1]}(M)$. From 
the long exact sequence associated to $K^\vee[2] \to \tau^{[-2,0]}(\HOM_{{\pr_M}}(E,E)[1])^\vee \to K$, we see 
the map $\tau^{[-2,0]}(\HOM_{\pr_M}(E,E)[1])^\vee \to K$ induces an isomorphism on $h^0$ and a surjection on $h^{-1}$.
As the obstruction theory $\phi_E$ factors through $K \to \bL_M$, we conclude
$K \to \bL_M$ is an obstruction theory.

Arguing as in the proof of Lemma \ref{lemmaori}, the identification of virtual classes follows from the construction
of $\otv{M}$ and \cite[Lemma 3.5, 3.6]{OT}.
\end{proof}

\begin{rmk}
The Proposition also follows from \cite[Proposition 1.18]{Par}, as the proof shows that
the assumptions on $K$ in loc.~cit. are satisfied.
\end{rmk}
 
 \subsection{Existence of maximal isotropic complexes}  \label{maxiso}
 \begin{ass} \label{lego}
Let $Y$ be a qcqs scheme,  $V,W \in \perf(Y)$ of tor-amplitude in $[-2,0]$. 
Assume there exists a non-degenerate 2-form $\vartheta\colon V^\vee[2] \simeq V$,
and $\sigma\colon V \to W$ is maximal isotropic, i.e., 
\begin{equation}
\label{ff}
W^\vee[2] \xrightarrow{\vartheta\circ \sigma^\vee[2]} V \xrightarrow{\sigma} W
\end{equation}
is a distinguished triangle. 
Further assume $h^{-2}(W)$ is locally free and $\tau^{\geqslant-1}W$ is perfect of tor-amplitude in $[-1,0]$.
\end{ass}
In this subsection, we will construct
another maximal isotropic complex perfect of tor-amplitude in $[-1,0]$ out of $\sigma$ under the assumptions above.

The triangle \eqref{ff} determines a map $W[-1] \to W^\vee[2]$, consider the diagram
\begin{equation} \label{hehe}
\xymatrix{
  & h^{-2}(W)[1] \ar[r]  
    & W[-1] \ar[r]\ar[d]  
      &  \tau^{\geqslant-1}W[-1]  \ar[r] 
         & h^{-2}(W)[2]\\
 h^{-2}(W)^\vee[-1]\ar[r]
   &(\tau^{\geqslant-1}W)^\vee[2]  \ar[r]    
      &   W^\vee[2] \ar[r]  
         &  h^{-2}(W)^\vee ,
}
\end{equation}
where the horizontal triangles come from shifting and dualizing
\[
h^{-2}(W)[2]   \to W \to \tau^{\geqslant-1}W.
\]

\begin{lem}
Let $Y$ be a qcqs scheme, and $W \in \perf(Y)$ of tor-amplitude in $[-2,0]$.
Assume $h^{-2}(W)$ is locally free and $\tau^{\geqslant-1}W$ is perfect of tor-amplitude in $[-1,0]$.
Then any map $W[-1] \to W^\vee[2]$ determines a unique map between triangles in
\eqref{hehe}.
\end{lem}

\begin{proof}
Apply $\hom_{\qcoh(Y)}(h^{-2}(W)[1], -)$ to the second row, as any map from $h^{-2}(W)[1]$ to $h^{-2}(W)^\vee$ or $h^{-2}(W)^\vee[-1]$ is zero, 
we have an isomorphism
\[
\hom_{\qcoh(Y)}(h^{-2}(W)[1], (\tau^{\geqslant-1}W)^\vee[2]) \simeq
\hom_{\qcoh(Y)}(h^{-2}(W)[1],  W^\vee[2]),
\]and there exists a unique map $h^{-2}(W)[1] \to (\tau^{\geqslant-1}W)^\vee[2]$ with which
the left square is commutative. 

We then have a map 
\[
\tau^{\geqslant-1}W[-1]  \to h^{-2}(W)^\vee
\] as part of the maps induced between the triangles. Its uniqueness can be seen by
applying $\hom_{\qcoh(Y)}(-,  h^{-2}(W)^\vee)$ to the first row.
\end{proof}
We define $L$ such that $$h^{-2}(W)[1] \to   (\tau^{\geqslant-1}W)^\vee[2] \to  L  \to h^{-2}(W)[2] $$ is a 
distinguished triangle. Applying the 3 by 3 lemma
\cite[Lemma 2.6]{may}
to the diagram below formed by the solid arrows, we obtain an object $K$ and dashed arrows including
$\xi$ and $\zeta$.
Note that the 3rd row and the 3rd column are distinguished triangles and they receive maps from existing distinguished triangles.
\[
\xymatrix{
h^{-2}(W)[1] \ar[r]\ar[d]    
  &  (\tau^{\geqslant-1}W)^\vee[2] \ar[d]\ar[r]           
    &   L \ar@{-->}[d]^-{\xi} \ar[r] 
      &  h^{-2}(W)[2] \ar[d]\\
W[-1]           \ar[r]    \ar[d] 
  & W^\vee[2] \ar[r]^-{\vartheta\circ \sigma^\vee[2]} \ar[d]
    &  V \ar@{-->}[d]^{\zeta} \ar[r]
	  & W \ar[d]\\
\tau^{\geqslant-1}W[-1] \ar[d]  \ar@{-->}[r]       
  &  h^{-2}(W)^\vee \ar@{-->}[r] \ar[d]
    & K \ar@{-->}[d] \ar@{-->}[r]
      & \tau^{\geqslant-1}W \ar[d]\\
h^{-2}(W)[2] \ar[r]
  & (\tau^{\geqslant-1}W)^\vee[3]\ar[r]
  	& L[1] \ar[r]
	  &h^{-2}(W)[3]    
}
\]

\begin{thm} \label{K}
Under Assumption \ref{lego}, construct $L$ and $K$ as above.
Then the map 
\[
V\xrightarrow{\xi^\vee[2]\circ \vartheta^{-1}} L^\vee[2]
\] is isomorphic to 
$V\xrightarrow{\zeta} K$, and $V \to K$ is maximal isotropic.

\end{thm}

\begin{proof}
Consider the triangle $L \xrightarrow{\xi} V \xrightarrow{\zeta} K$ and 
its dual with shift
 $$K^\vee[2] \xrightarrow{ \vartheta \circ \zeta^\vee[2]} V 
 \xrightarrow{\xi^\vee[2]\circ \vartheta^{-1}} L^\vee[2].$$
We first show 
the compositions $L \xrightarrow{\xi} V \xrightarrow{\xi^\vee[2]\circ \vartheta^{-1}} L^\vee[2]$ 
and $K^\vee[2] \to V \to K$ are zero.

The cohomology sheaves of $K$ and $L$ are related to those of $V$ as follows.
In the diagram
\[
\xymatrix{
  (\tau^{\geqslant-1}W)^\vee[2] \ar[d]\ar[r]         
    &   L \ar@{-->}[d]^-{\xi} \ar[r] 
      &  h^{-2}(W)[2] \ar[d]\\
W^\vee[2] \ar[r] \ar[d]
    &  V \ar@{-->}[d]^{\zeta} \ar[r]
	& W \ar[d]\\
 h^{-2}(W)^\vee \ar@{-->}[r] 
    & K  \ar@{-->}[r]
      & \tau^{\geqslant-1}W,
}
\] 
All rows and columns are distinguished triangles. From the first row and 3rd row, we see $L$ is perfect of tor-amplitude in $[-2,-1]$,
$K$ perfect of tor-amplitude in $[-1,0]$.
From the induced long exact sequence on cohomology  of the middle column, we get isomorphisms
$h^{-2}(L) \simeq h^{-2}(V)$, $h^0(V)\simeq h^0(K)$ and a short exact sequence 
\[
0 \to h^{-1}(L) \to h^{-1}(V) \to h^{-1}(K) \to 0.
\]

Now we dualize the diagram. Apply $(-)^\vee[2]$ to the diagram and identify $V^\vee[2]$ with $V$ using $\vartheta$, 
we get another diagram
\[
\xymatrix{
  (\tau^{\geqslant-1}W)^\vee[2] \ar[d]     \ar@{-->}[r]  
    & K^\vee[2] \ar@{-->}[r] \ar@{-->}[d]^-{\vartheta\circ \zeta^\vee[2]}
    		&  h^{-2}(W)[2] \ar[d]\\
W^\vee[2] \ar[r] \ar[d]
    &  V\ar[r]^-{-\sigma} \ar@{-->}[d]^-{\xi^\vee[2]\circ \vartheta^{-1}}
		& W \ar[d]\\
 h^{-2}(W)^\vee  \ar[r]
    & L^\vee[2]\ar[r]
      & \tau^{\geqslant-1}W. 
}
\]
The negative sign comes from $\vartheta^\vee[2]=-\vartheta$, and can be absorbed into $W$ if we
further apply the isomorphism $-\id$ to the the right column.
Then both diagrams by construction contain
\[
\xymatrix{
  (\tau^{\geqslant-1}W)^\vee[2] \ar[d]        
    & 
      &  h^{-2}(W)[2] \ar[d]\\
W^\vee[2] \ar[r] \ar[d]
    &  V \ar[r]
	& W \ar[d]\\
 h^{-2}(W)^\vee  
    & 
      & \tau^{\geqslant-1}W.
}
\] 
Consider the composition $L \xrightarrow{\xi} V \xrightarrow{\xi^\vee[2]\circ \vartheta^{-1}} L^\vee[2]$, it is sandwiched between a map of triangles
\[
\xymatrix{
  (\tau^{\geqslant-1}W)^\vee[2] \ar[d]^-{0}  \ar[r]   
    & L\ar[r] \ar@{-->}[d]
      &  h^{-2}(W)[2] \ar[d]^-{0}\\
 h^{-2}(W)^\vee \ar[r] 
    & L^\vee[2]\ar[r]
      & \tau^{\geqslant-1}W.
}
\]
Using the short exact sequence
\[
\hom(L, h^{-2}(W)^\vee)=0 \to \hom(L, L^\vee[2]) \to \hom(L, \tau^{\geqslant-1}W),
\]
it is easy to conclude 
$L \to L^\vee[2]$ is zero. Same argument shows $K^\vee[2] \to V \to K$ is zero.

As $L \to V \to L^\vee[2]$ is zero, and $L \to V \to K$ is a distinguished triangle,
we obtain a map $K \to L^\vee[2]$ under $V$. Similarly using $K^\vee[2] \to V \to K$ being zero, we obtain a map
$L^\vee[2] \to K$ under $V$. The long exact sequence on cohomology sheaves then shows 
the maps $K \to L^\vee[2]$ and $L^\vee[2] \to K$ induces isomorphism on $h^{-1}$ and $h^0$, because 
 $h^{i}(V) \to h^{i}(L^\vee[2]), h^{i}(V) \to h^{i}(K),i=-1,0$ are all surjective. As $K$ and $L^\vee[2]$
 are perfect in $[-1,0]$, we see $V \to K$ is isomorphic to $V \to L^\vee[2]$.
\end{proof}

\section{Applications to Calabi-Yau 4-folds fibered over curves}

\subsection{Virtual classes on Calabi-Yau 4-folds fibered over curves}

\begin{set} \label{overcurve}

Let $X$ be a (connected) smooth complex projective CY 4-fold and $\pi \colon X \to C$ a projective surjective map to a smooth projective curve $C$ with connected fibers.
Let $\beta\in H_2(X,\bZ)$ satisfy $\pi_*\beta=0$.
 Let $U \subset C$ be an open subset over which $\pi$ is smooth.
\end{set}

Consider the moduli spaces $M_\beta(\pi)$ and $M_\beta(X)$, which are isomorphic under (ref.~Proposition~\ref{relabs}): 
$$\epsilon\colon M_\beta(\pi) \to M_\beta(X). $$
%(reion \ref{moduli}).  
Denote $M$ to be the preimage
of $U$ under $M_\beta(\pi) \to C$,~i.e.~we have the following Cartesian diagram 
\begin{align*}
\xymatrix{
M \ar[r]^{ } \ar[d]  \ar@{}[dr]%|{\Box}  
	& M_\beta(X) \ar[d] \\
U \ar[r]^{ } 
	&  C.}
%\xymatrix{
% M \ar[r] \ar[d] & M_\beta(X) \ar[d] \\
%  U \ar[r] & C}
\end{align*}
Then $M$ is a moduli space of sheaves on the fibers of 
the smooth fibration $\pi^{-1}U \to U$. 
As $X$ is a smooth projective
CY 4-fold,  for any orientation on $M_\beta(X)$, we have $\otv{M_\beta(X)}$. 
For the smooth family of CY 3-folds
$\pi^{-1}U\to U$, we have 
the family of moduli spaces $M\to U$ and its 
Behrend-Fantechi virtual class $\bfv{M}$.
Identify $M_\beta(\pi)$ and $M_\beta(X)$ via $\epsilon$ and view $M$ as an open subscheme of $M_\beta(X)$
with imbedding $\iota\colon M \to M_\beta(X)$,
we will relate the two virtual classes.

Let $F$ be the normalized universal sheaf (Remark \ref{norm}) on $M\times_CX$,
then $F\in \perf(M\times_C X)$.
 As in Section \ref{functoriality}, we have the maps $j, p, q$ and sheaf $G=j_*F$. 
 %The family $G$ determines the open embedding $\iota\colon M \to M_\beta(X)$.
 Denote $r$ the projection 
$M_\beta(X) \times X \to M_\beta(X)$, and $\kappa\colon M\times X \to M_\beta(X) \times X$
the map $\iota\times \id_X$.
To summarize, we have a commutative diagram
\[
\xymatrix{
M\times_C X \ar[r]^-{j} \ar[dr]_-{p} 
	& M\times X \ar[r]^-{\kappa} \ar[d]^-{q} 
		& M_\beta(X)\times X\ar[d]^-{r}\\
{}  	
	&M \ar[r]^-{\iota}
		&M_\beta(X)	
}
\]

%\subsection{}

Denote $E$ the normalized universal sheaf on $M_\beta(X) \times X$. As $\kappa^*E=G=j_*F$, 
we have a commutative diagram
\begin{equation} \label{obscom}
\xymatrix{
\iota^* (\HOM_r(E,E)[1])^\vee \ar[d]^-{\iota^*\phi_E} \ar[r]^-{\sim} 
	&  (\HOM_q(G,G)[1])^\vee \ar[d]^-{\phi_G}	\ar[r]
		& (\HOM_p(F,F)[1])^\vee \ar[d]^-{\phi_F}\\
\iota^*\bL_{M_\beta(X)} \ar[r]^-{\sim} 
	&\bL_M\ar[r]
		& \bL_{M/C}	.
}
\end{equation}
Viewing $M$ as an open subset of $M_\beta(X)$, clearly the left square commutes\,\footnote{More generally this follows from functoriality for $\phi$ under pullback, see Proposition \ref{phipb}.}.
The right square is commutative by Proposition \ref{phipf}. The maps $\phi_F$ and $\phi_E$ \eqref{phie} are used to produce obstruction theories for $M$ and 
$M_\beta(X)$ respectively.  The diagram will be used to 
relate $\otv{M_\beta(X)}$ and $\bfv{M}$.

\begin{prop}\label{prop on can ori}
In Setting \ref{overcurve}, 
 the obstruction theory
 for $M_\beta(X)$
\[
\phi_E\colon \tau^{[-2,0]}(\HOM_r(E,E)[1])^\vee \to \bL_{M_\beta(X)}
\] has a canonical orientation on the open subscheme $M$.
\end{prop}

\begin{proof}
Note that the restriction of $\phi_E$ to $M$ is given by 
\[
\phi_G \colon \tau^{[-2,0]}(\HOM_q(G,G)[1])^\vee  \to \bL_M.\]
The canonical orientation will come from a maximal isotropic complex of
tor-amplitude in $[-1,0]$.

As $C$ is a smooth curve, $j\colon M \times_C X \to M \times X$ is the inclusion of an effective
Cartier divisor.
Therefore  we have a self-dual triangle
\[
j_*\HOM(F,F) \to \HOM(G,G) \to j_*\HOM(F,F\otimes \mathcal{N}[-1])
\] on $M\times X$ by Proposition \ref{selfdual}, where $\cN$ denotes the normal sheaf of $j$.
As the dualizing sheaf $\omega_X$ is trivial, the dualizing complex of $q\colon M \times X \to M$ is $\cO[4]$. Pushing forward along $q$ the self-dual triangle
and taking dual and shift we obtain
\begin{equation}\label{ori}
\HOM_p(F,F)[3] 
\to \HOM_q(G,G)[3] \simeq (\HOM_q(G,G)[1])^\vee
\to (\HOM_p(F,F)[1])^\vee.
\end{equation}

Let $V=\tau^{[-2,0]}(\HOM_q(G,G)[1])^\vee$ and $W=\tau^{[-2,0]}(\HOM_p(F,F)[1])^\vee$, 
then they are perfect of tor-amplitude in $[-2,0]$.
From \eqref{ori} we obtain by truncation and duality 
a triangle 
\[
W^\vee[2] \to V \to W,
\] and we see $V \to W$ is maximal isotropic.
As we are in the situation of Assumption \ref{lego},
the proposition follows from Theorem \ref{K} and Lemma \ref{lemmaori}.
Note that 
\[
\tau^{\geqslant -1}W=\tau^{[-1,0]}(\HOM_p(F,F)[1])^\vee
\] is perfect of tor-amplitude
 $[-1,0]$ (see e.g., \cite[4.4]{HT}), and $h^{-2}(W)$ is locally free by Lemma \ref{lem on wc}.
\end{proof}

\begin{rmk}
Let $\bF$ be the normalized universal sheaf on $M_\beta(\pi) \times_C X$, it is not clear to us
whether the derived pushforward of $\HOM(\bF,\bF)$ to $M_\beta(\pi)$ is perfect if $\pi$ has singular fibers.
\end{rmk}

%When $F$ is stable, the natural maps $\cO_M \to p_*\HOM(F,F)$  induces isomorphisms on $h^0$, and its cones
%$\tau^{>0}p_*\HOM(F,F)$ is perfect of tor-amplitude in $[1,3]$.
%Taking dual, we see 
%$(p_*\HOM(F,F))^\vee\to \cO_M$ induces an isomorphism on $h^0$ and 
%\[
%(\tau^{>0}p_*\HOM(F,F))^\vee \simeq \tau^{<0}((p_*\HOM(F,F)^\vee).
%\]
%Similarly, 
%$\cO_M \to q_*\HOM(G,G)$ induces an isomorphism on $h^0$,
%and $\tau^{>0}p_*\HOM(G,G)$ is perfect of tor-amplitude in $[1,4]$. Dually we obtain 
%$(q_*\HOM(G,G))^\vee \to \cO_M$, and 
%$(\tau^{>0}q_*\HOM(G,G))^\vee \simeq \tau^{<0}((q_*\HOM(G,G))^\vee)$.

\begin{lemma}\label{lem on wc}
Denote
$W=\tau^{[-2,0]}(\HOM_p(F,F)[1])^\vee$, then $h^{-2}(W)\simeq {\Omega_C^\vee}_{|_M}$.
It follows that $h^0(W^\vee[2])\simeq {\Omega_C}_{|_M}$. Here ${\Omega_C}_{|_M}$ denotes
the pullback of $\Omega_C$ to $M$.
\end{lemma}

\begin{proof}

By duality, we have 
\[
(p_*\HOM(F,F)[1]))^\vee \simeq p_*\HOM(\HOM(F, F)[1], \omega_p^\bullet)
\simeq \HOM_p(F, F\otimes \omega_p^\bullet[-1]).
\]
From the Cartesian diagram
\[
\xymatrix{
M\times_CX \ar[r]^-{\pr_X}\ar[d]^-{p} 
	& X \ar[d]^-{\pi}\\
M\ar[r]
	&C
}
\] where $\pi$ is flat, we see $\omega_p^\bullet =p^!\cO_M\simeq  \pr_X^*\pi^!\cO_C =\pr_X^* \omega_\pi^\bullet$
(\cite[\href{https://stacks.math.columbia.edu/tag/0AA8}{0AA8},
\href{https://stacks.math.columbia.edu/tag/0A9P}{0A9P}]{stacks}).

As $\omega_X^\bullet \simeq \cO[4]$, and 
$\omega_X^\bullet \simeq \omega_\pi^\bullet \otimes \omega_C^\bullet$,
we see $\omega_\pi^\bullet \simeq \Omega_C^\vee[3]_{|_X}$.
By the projection formula,
\[
\HOM_p(F, F\otimes \omega_p^\bullet[-1]) \simeq
\HOM_p(F, F)\otimes {\Omega_C}^\vee_{|_M}[2],
\] and we see 
$h^{-2}(W)\simeq {\Omega_C^\vee}_{|_M}$, taking dual gives 
$h^0(W^\vee[2])\simeq {\Omega_C}_{|_M}$.
\end{proof}

\begin{thm} \label{OTBF}
In Setting \ref{overcurve},  consider moduli spaces $M, M_\beta(X)$ and their virtual classes as above.
Assume there exists an orientation of $M_\beta(X)$ that restricts to the canonical orientation
of $M$, then with that orientation, we have
\[
\iota^*\otv{M_\beta(X)} = \bfv{M},
\]
where $\iota\colon M \to M_\beta(X)$ denotes the open immersion.  
\end{thm}

\begin{proof}

Let $V=\tau^{[-2,0]}(\HOM_q(G,G)[1])^\vee$ and $W=\tau^{[-2,0]}(\HOM_p(F,F)[1])^\vee$.
Consider the commutative diagram between horizontal distinguished triangles
\[
\xymatrix{
  (\tau^{\geqslant-1}W)^\vee[2] \ar[d]\ar[r]         
    &   L \ar[d]^-{\xi} \ar[r] 
      &  h^{-2}(W)[2] \ar[d]\\
W^\vee[2]\ar[r]\ar[d]
	&V  \ar[d]^-{\phi_G} \ar[r]
		& W \ar[d]^-{\phi_F}\\
{\Omega_C}_{|_M}\ar[r]
	&	\bL_M^{\geqslant-1}\ar[r]
		& \bL_{M/C}^{\geqslant-1}
}
\]
Here the top half is obtained by the argument constructing $L$ in Section \ref{maxiso},
and the lower half is induced by the lower right commutative square obtained by truncating
\eqref{obscom}.
By Theorem \ref{K}, there exists 
a distinguished triangle $L \to V \to K$ and 
$V \to K$ is maximal isotropic.

As the cone of $\phi_G,\phi_F$
lie in $D^{\leqslant -2}(M)$, the
cone of the non-unique map $W^\vee[2] \to{\Omega_C}_{|_M}$
belongs to $D^{\leqslant -1}(M)$, and $W^\vee[2] \to {\Omega_C}_{|_M}$ is surjective on $h^0$.
As $h^0(W^\vee[2])\simeq {\Omega_C}_{|_M}$ by Lemma \ref{lem on wc}, the surjection is an isomorphism, and we conclude the left vertical triangle is distinguished.

Now it is easy to see the composition of  $L \to V \to \bL_M^{\geqslant-1}$ is zero, as there are no nonzero maps from $h^{-2}(W)[2]$ to $\bL_M^{\geqslant -1}$, and the composition of $(\tau^{\geqslant-1}W)^\vee[2] \to L$
and $L  \to \bL_M^{\geqslant-1}$ factors through the left column.
Then $V \to \bL_M^{\geqslant-1}$ factors through the maximal isotropic $V \to K$,
and $\iota^*\otv{M_\beta(X)}$ is the class defined using the perfect obstruction theory
$K \to \bL_M^{\geqslant-1}$ by Proposition \ref{OTviaBF}.

Note that $\tau^{[-1,0]}W \to \bL_{M/C}^{\geqslant-1}$
is the relative perfect obstruction theory used to define $\bfv{M}$. 
As in \cite[Appendix B]{GP}, we obtain an absolute perfect obstruction theory 
 $K' \to \bL_M^{\geqslant-1}$ and a 
diagram
\[
\xymatrix{
{\Omega_C}_{|_M}\ar[r]\ar[d]
	& K'\ar[r]\ar[d]
		& \tau^{[-1,0]}W \ar[d]^-{\phi_F}\\
{\Omega_C}_{|_M}\ar[r]
	&	\bL_M^{\geqslant-1}\ar[r]
		& \bL_{M/C}^{\geqslant-1}.
}   
\]  The obstruction theory $K' \to \bL_M^{\geqslant-1}$ produces the class $\bfv{M}$.

By the construction of $K$, we have a distinguished triangle
$
{\Omega_C}_{|_M} \to K \to \tau^{[-1,0]}W
$, and we see the class of $K$
and $K'$ in $K_0(\perf(M))$ are the same.
It then follows from Siebert's formula \cite{Sie} that the classes defined using
$K \to \bL_M^{\geqslant-1}$ and $K' \to \bL_M^{\geqslant-1}$ are the same.
\end{proof}

\begin{rmk}
The theorem is a generalization of \cite[Lemma 2.6]{CMT1}. In loc.~cit.,
$X=Y\times E$ is the product of a CY 3-fold $Y$ and an elliptic curve $E$, and $\pi$ the projection
to $E$, then $\pi$ is smooth and we have $M=M_\beta(X)$.
\end{rmk}

\begin{rmk}
Without the assumption on orientation as in \cite[Theorem 3.1]{Par},  
$\iota^*\otv{M}$ and $\bfv{M}$ agree up to a $\pm$ sign on each connected component of $M$.
\end{rmk}

\begin{rmk}
If the fibers of $M \to U$ are all connected, then $M$ is connected
(\cite[\href{https://stacks.math.columbia.edu/tag/0377}{0377}]{stacks})
and the assumption on orientation holds.  Elliptic fibrations
constructed as Weierstrass models are such examples (see e.g., \cite[Section 2.1]{CMT1} and take $Y$ the product of  $\bP^1$ and a Del Pezzo surface instead of $\bP^3$).
%Let $Y=\bP^1 \times S$ be a smooth connected Fano 3-fold, 
%and $X \to Y$ an elliptic fibration constructed as a Weierstrass model so that the singualr fibers of $X \to Y$ are integral curves of arithmetic genus 1 with a single node or cusp.  
%For the fibration $X \to Y \xrightarrow{\pr_{\bP^1}} \bP^1$, and $\beta$ a multiple of the class of 
%an elliptic fiber.
\end{rmk}

 Let $c \in C$ be a closed point. Denote  $X_c$ the fiber of $\pi$ over $c$ with inclusion $i: X_c\to X$, $M_\beta(X_c)$ the fiber of $M_\beta(\pi) \to C$ over $c$.
Note that 
$$M_\beta(X_c)=\bigsqcup_{\begin{subarray}{c}\beta' \in H_2(X_c, \mathbb{Z}) \\ i_*\beta'=\beta\end{subarray}} M_{\beta'}(X_c)$$ 
is the disjoint union of $M_{\beta'}(X_c)$ over all $\beta'\in H_2(X_c)$
whose image in $H_2(X)$ is $\beta$. If $c\in U$, then $X_c$ is a smooth CY 3-fold and the inclusion $i_c\colon c \to C$ is a regular closed imbedding, in which case denote $\bfv{M_\beta(X_c)}$ the virtual class of $M_\beta(X_c)$.

For $c\in U$, one can pull back the
complex $\tau^{[-2,0]}(\HOM_q(G,G)[1])^\vee$
with canonical orientation in Proposition \ref{prop on can ori} to 
$M_\beta(X_c)$, which we call the \textit{canonical orientation of} $M_\beta(X_c)$. 
 \begin{cor} \label{DT43}
 Let $c\in U$. 
Consider the Cartesian diagram
\[
\xymatrix{
M_\beta(X_c) \ar[r]\ar[d] 
	& M_\beta(\pi)\ar[d]\\
c \ar[r]^-{i_c} 
	& C.
} 
\]
Assume there exists an orientation of $M_\beta(X)$ that restricts to the canonical orientation
of $M_\beta(X_c)$. Then with that orientation, we have 
%If different connected components of $M_\beta(X_c)$ are in different connected components of $M_\beta(X)$, 
%then there exists an orientation on $M_\beta(X)$ 
% $M_\beta(X_c)$ is connected, 
\[
i_c^!\otv{M_\beta(X)}= \bfv{M_\beta(X_c)}.
\] Here $i_c^!$ denotes the refined Gysin pullback along the regular closed immersion $i_c$.
\end{cor}

\begin{proof}
Factor $i_c$ as $c\to U \to C$, and the corollary follows from the proof of Theorem \ref{OTBF} and deformation invariance of virtual classes for the family $M \to U$ (\cite[Cor.~4.3]{HT}). 
Note that only the connected component(s) of $M$ containing $M_\beta(X_c)$ are relevant for the argument.
\end{proof}
\begin{rmk}
If different connected components of $M_\beta(X_c)$ are in different connected components of $M_\beta(X)$, then 
the assumption of Corollary \ref{DT43} is satisfied. 
\end{rmk}

\subsection{Numerical results}

We turn to numerical consequences of Corollary \ref{DT43}.
Consider the Cartesian diagram
\[
\xymatrix{
M_\beta(X) \times X \ar[r]^-{\pr_X} \ar[d]^-{q} 
	& X\ar[d] \\
M_\beta(X) \ar[r]                         
	& {*}.
}
\]
Using the normalized universal sheaf $\bG$ on $M_\beta(X) \times X$ as a Fourier-Mukai kernel, we
have tautological complexes and descendants as follow.

For any $L\in \Pic(X)$, we define its tautological complex
\begin{align}\label{def of tau cpx}
L^{[1]}:= \mathbf{R}q_*(\bG\otimes \pr_X^* L)\in \perf(M_\beta(X)).
\end{align}
%Denote $n_{0,\beta}(X,L)=\int_{\otv{M_\beta(X)}}c_1(L^{[1]})$.
%\footnote{$\pr_X\colon M_beta(X) \times_BX$}
\begin{rmk}
Denote $p\colon M_\beta(X) \times_C X\to M_\beta(X)$ the projection map,
and $\bF$ the normalized universal sheaf on $M_\beta(X) \times_C X$.
It follows from the definition
that 
\[
L^{[1]}=\mathbf{R}p_*(\bF\otimes \pr_X^* L).
\]
%If $L=\pi^*\cO_C(c)$ for some $c\in C$, by projection formula, we have $$L^{[1]}=\mathbf{R}q_*$$
\end{rmk}
Let $k,m$ be non-negative integers and  $\gamma \in H^{m}(X,\bQ)$ a cohomology class,
the map  
 \[
 q_*(\pr_X^* \gamma \cap \ch_{3+k}(\bG)\cap q^*(-))\colon 
H_l(M_\beta(X),\bQ) \to H_{l+2-2k-m}(M_\beta(X),\bQ)
\] 
determines 
\begin{equation}\label{equ on tau}\tau_k(\gamma) \in H^{2k+m-2}(M_\beta(X),\bQ)\end{equation} 
so that the map
is given by capping with $\tau_k(\gamma)$.
Denote
\[
\langle \tau_k(\gamma)\rangle_\beta =\deg (\tau_k(\gamma)\cap \otv{M_\beta(X)}).
\footnote{we use the same notation $\otv{M_\beta(X)}$ for its class in $H_2(M_\beta(X),\bQ)$}
\]

\begin{rmk}
As $\otv{M_\beta(X)} \in H_2(M_\beta(X))$, $\langle \tau_k(\gamma)\rangle_\beta$
is zero unless $2k+m=4$ and $0\leqslant k\leqslant2$.
\end{rmk}

\begin{rmk}
If $\gamma$ belongs to the Chow ring of $X$, $\tau_k(\gamma)$ is in the operational Chow ring of $M_\beta(X)$.
\end{rmk}

\begin{lemma}\label{lem on tau}
We have the following: 
\begin{enumerate}
	\item Let $L=\pi^*\cO_C(c)$, $c\in C$ and $\gamma \in H^m(X,\bQ)$, then 
		\[
		\tau_k(\gamma \cup c_1(L)) = \tau_k(\gamma)\cap c_1(\cO_C(c)_{|_{M_\beta(X)}}).
		\] Here $\cO_C(c)_{|_{M_\beta(X)}}$ denotes the pullback of $\cO_C(c)$ to $M_\beta(\pi) \simeq M_\beta(X)$ under the 
		map	$M_\beta(\pi)\to C$.
	\item We have $\tau_1(1_X)=1$, where $1_X\in H^*(X)$ is the identity, and
		$\tau_0(\gamma)=\int_\beta \gamma $ for $\gamma\in H^2(X)$.
		\end{enumerate}
\end{lemma}

\begin{proof}
(1)
Because the universal sheaf $\bG$ is supported in $M_\beta(X) \times_C X$, the Chern character $\ch(\bG)$ comes from a localized Chern character
that maps a homology class of $M_\beta(X) \times X$ to one of $M_\beta(X) \times_C X$ (see e.g., \cite{Iv}).
Then the identification follows since the restriction of $L$ to $M_\beta(X) \times_C X$
is $q^*\cO_C(c)_{|_{M_\beta(X)}}$.

(2) 
Let $\pt$ be a closed point of $M_\beta(X)$ determined by a stable sheaf $\cF$ on $X$, and
$[\pt]$ its class in $H_0(M_\beta(X))$.
Note that  $\tau_1(1_X)$ and $\tau_0(\gamma),\gamma\in H^2(X)$ are of degree zero
and acts on $[\pt]$ as scalar multiplication by $\int_X \ch_4(\cF)$ and $\int_X\ch_3(\cF) \cup \gamma$. 
The results  then follow from  $\int_X \ch_4(\cF)=\chi(\cF)=1$ and $\int_X\ch_3(\cF) \cup \gamma=\int_\beta \gamma$.
\end{proof}
\begin{thm} \label{222}
Let $c\in U$, $L=\pi^*\cO_C(c)$ and $\gamma \in H^2(X, \bQ)$. Assume
\[
i_c^!\otv{M_\beta(X)}= \bfv{M_\beta(X_c)},
\] 
which is satisfied if there exists an orientation of $M_\beta(X)$ that restricts to the canonical orientation
of $M_\beta(X_c)$ by Corollary \ref{DT43}. 
Then we have the following identities:
$$\langle\tau_0(\gamma \cup c_1(L) )\rangle_{\beta}=
\int_\beta  \gamma\cdot \int_{[M_\beta(X_c)]^{\mathrm{vir}}}1
=\int_\beta  \gamma\cdot \langle\tau_1( c_1(L) )\rangle_{\beta},
$$
\[\frac{1}{12}\langle\tau_0(c_2(X)) \rangle_{\beta}+\langle\tau_2(1_X) \rangle_{\beta}=0.
\]
%\[\langle\tau_0(\gamma \cup c_1(L) )\rangle_{\beta}=\int_\beta  \gamma \cdot \deg \bfv{M_\beta(X_c)},\]
\end{thm}

\begin{proof}
By Lemma \ref{lem on tau}, Corollary \ref{DT43} and  
\[
c_1(\cO_C(c)_{|_{M_\beta(X)}}) \cap \otv{M_\beta(X)} 
=\deg i_c^!\otv{M_\beta(X)},
\]
we obtain 
\begin{equation}\label{pf equ1}\langle\tau_0(\gamma \cup c_1(L) )\rangle_{\beta}=\int_\beta  \gamma\cdot \int_{[M_\beta(X_c)]^{\mathrm{vir}}}1.
\end{equation}
%For $L=\pi^*\cO_C(c)$, the second equation comes from computating 
%$ \int_{\otv{M_\beta(X)}}c_1(L^{[1]}) $  in the Chow ring in two ways .
By a direct calculation, we have $$\det(L^{[1]})=\cO_C(c)_{|_{M_\beta(X)}}, $$ 
and hence 
\begin{equation}\label{pf equ2}
\int_{\otv{M_\beta(X)}}c_1(L^{[1]})=\int_{\otv{M_\beta(X)}}c_1(\det(L^{[1])})=\deg \bfv{M_\beta(X_c)}.
\end{equation}
Apply the Grothendieck Riemann-Roch theorem \cite[Theorem 18.2 (3)]{Ful} to the smooth map
$q\colon M_\beta(X) \times X \to M_\beta(X)$, 
we have
\begin{align}\label{pf equ3}
\int_{\otv{M_\beta(X)}}c_1(L^{[1]})
&=\int_{\otv{M_\beta(X)}}\ch(L^{[1]}) \\ \nonumber 
&=\int_{\otv{M_\beta(X)}} q_*\left(\ch(\bG) \cdot \pr_X^*(\ch(L)\cdot \td(X))\right) \\ \nonumber
&=\int_{\otv{M_\beta(X)}} q_*\left(\ch_3(\bG)\cdot \pr_X^*\left(\frac{1}{2}c^2_1(L)
   +\frac{1}{12}c_2(X)\right)\right) \\ \nonumber
&\, +\int_{\otv{M_\beta(X)}} q_*\left(\ch_4(\bG)\cdot \pr_X^*\left(c_1(L)\right)\right)
       +\int_{\otv{M_\beta(X)}} q_*\left(\ch_5(\bG)\right)   \\ \nonumber
&=\frac{1}{2}\langle\tau_0(c^2_1(L)) \rangle_{\beta}+\frac{1}{12}\langle\tau_0(c_2(X)) \rangle_{\beta}+\langle\tau_1(c_1(L)) \rangle_{\beta}+\langle\tau_2(1_X)\rangle_{\beta}.
\end{align}
Since $L=\pi^*\cO_C(c)$, then $c_1^2(L)=\pi^*c_1^2(\cO_C(c))=0$,
while 
\[
\tau_1(c_1(L))=\tau_1(1_X)\cap c_1(\cO_C(c)_{|_{M_\beta(X)}})=c_1(\cO_C(c)_{|_{M_\beta(X)}}).
\]

It follows that 
\begin{equation}\label{pf equ4}\frac{1}{2} \langle\tau_0(c^2_1(L)) \rangle_{\beta}=0, \quad \langle\tau_1(c_1(L)) \rangle_{\beta}=\deg \bfv{M_\beta(X_c)}, \end{equation}

Combining
\eqref{pf equ2}, \eqref{pf equ3},
and \eqref{pf equ4}, we see
\begin{equation}\label{pf equ5}
%\frac{1}{2}\langle\tau_0(c^2_1(L)) \rangle_{\beta}+
\frac{1}{12}\langle\tau_0(c_2(X)) \rangle_{\beta} + \langle\tau_2(1_X)\rangle_{\beta}=0.  \end{equation}

\end{proof}
%\begin{rmk}
%The first equation verifies \cite[Conjecture 2.4]{CMT1},
%and the second equation was conjectured in \cite[\S 1.7]{CT2}.
%\end{rmk}

\subsection*{Acknowledgments}
F.Q. thanks Jeongseok Oh and Weizhe Zheng for correspondence.
We thank the referees for comments and suggestions that lead to improved exposition.
Y.C.~is partially supported by RIKEN Interdisciplinary Theoretical and Mathematical Sciences
Program (iTHEMS), JSPS KAKENHI Grant Number JP19K23397 and Royal Society Newton International Fellowships Alumni 2022 and 2023. 
F.Q. is partially supported by NSFC grant 11801185.
%\subsection*{Statements and Declarations}
%We have no conflicts of interest to disclose.

\begin{appendix}

\section{}
Let $i\colon Z \to X$ be the inclusion of an effective Cartier divisor into an arbitrary scheme $X$, and $\cN$ the normal
bundle of $i$. This section concerns the distinguished triangle \eqref{dist} constructed
from adjunctions involving $i_*$.

\begin{lem}\label{01}
Let $F \in \rmd(Z)$. We have
\[
i_*i^*i_*F \simeq i_*(F\otimes \cN^\vee)[1]\oplus i_*F.
\]
It follows that when $\cF$ is an $\cO_Z$-module,
$h^l(\mathbf{L}i^*i_*\cF)=0$ unless $l=-1,0$, and $h^0(\mathbf{L}i^*i_*\cF)=\cF$,
$h^{-1}(\mathbf{L}i^*i_*\cF)=\cF\otimes \cN^\vee$.
\end{lem}

\begin{proof}
The proof is the same as \cite[Corollary 11.4 (iii)]{Huy} for $F \in \coh^b(Z)$.
 For the closed immersion $i$,
the projection formula 
$
i_*F\otimes E \simeq i_*(F\otimes i^*E)
$ holds for any $E \in \rmd(X)$ and $F\in \rmd(Z)$ (\cite[\href{https://stacks.math.columbia.edu/tag/0B55}{0B55}]{stacks}). Then for $F\in \rmd(Z)$, we have
\[i_*i^*i_*F
\simeq i_*(i^*i_*F \otimes \cO_Z)
\simeq  i_*F\otimes i_*\cO_Z 
 \simeq 
 i_* (F \otimes 
\mathbf{L}i^*i_*\cO_Z)
\] and $\mathbf{L}i^*i_*\cO_Z\simeq i^*\{\cO_X(-Z) \to \cO_X\} \simeq \cN^\vee[1]\oplus\cO_Z$.
\end{proof}

\begin{lem} \label{triangle}
If $\cF$ is a finite type quasi-coherent sheaf, then the map $\cF\otimes \cN^\vee[1]\to \mathbf{L}i^*i_*\cF$ induces an isomorphism on $h^{-1}$, and $\mathbf{L}i^*i_*\cF \to \cF$ induces an isomorphism on $h^0$.
It follows from Lemma \ref{01} that 
\[\cF\otimes \cN^\vee[1]\to \mathbf{L}i^*i_*\cF \to \cF\]
is isomorphic to the distinguished triangle
\[
\tau^{\le -1}\mathbf{L}i^*i_*\cF
\to
\mathbf{L}i^*i_*\cF
\to
\tau^{\ge 0}
\mathbf{L}i^*i_*\cF.
\] %T
\end{lem}

\begin{proof}

Adjunctions give rise to retractions 
\[
i_*\cF \to i_*i^!i_*\cF \to i_*\cF
\]
and 
\[
i_* \cF\to i_*\mathbf{L}i^*i_*\cF \to i_*\cF.
\]
Taking $h^{0}$
we obtain retractions $\cF \to h^{0}(i^!i_*\cF) \to \cF$ and $\cF \to h^0(\mathbf{L}i^*i_*\cF) \to \cF$.
As 
\[
h^{0}(i^!i_*\cF)\simeq h^{0}(\mathbf{L}i^*i_*\cF\otimes \cN[-1])\simeq 
h^{-1}(\mathbf{L}i^*i_*\cF)\otimes \cN\simeq\cF\]
and $h^0(\mathbf{L}i^*i_*\cF)\simeq \cF$, we have two retractions of $\cF$ into itself. Locally,
$\cF$ is given by a finitely generated module, it then follows from Nakayama's lemma that
such retractions must be isomorphisms.
\end{proof}

\begin{rmk} \label{qcs}
Assume $X$ is quasi-compact and separated, then  
we have equivalence between the categories of Fourier-Mukai kernels and continuous maps between 
stable $\infty$ or dg enhancement of $\qcoh(Z)$ and $\qcoh(X)$
(\cite[Theorem 1.2 (2)]{BZFN} \cite[Corollary 1.8]{To}).
For $F \in \qcoh(Z)$, 
to show adjunctions induce a distinguished triangle 
\[F\otimes \cN^\vee[1]\to i^*i_*F \to F,\]
we need to show that, for the corresponding maps between Fourier-Mukai kernels,
\[
{\Delta}_*\cN^\vee[1] \to \{ {\cI_X}_{|_{Z\times Z}} \to \cO_{Z\times Z} \}
\to {\Delta}_*\cO_Z
\] form a distinguished triangle.
It is not immediately clear that this sequence can be identified with Eq.\eqref{FM}.
From Eq.\eqref{FM}, we see the complex
${\cI_X}_{|_{Z\times Z}} \to \cO_{Z\times Z}$
only has none trivial cohomology in degree $-1$ and $0$, 
and ${\Delta}_*\cN^\vee[1] $
and ${\Delta}_*\cO_Z$ are its cohomology sheaves in degree $-1$ and $0$ respectively.
To make the identification, 
argue as in the proof of the
previous lemma using the retractions $i_* \to i_*i^!i_* \to i_*$ and $i_* \to i_*i^*i_* \to i_*$.

\end{rmk}

\section{}

In this section, several identification results used in Sections \ref{sect on pre} \& \ref{sect on pre2}
are elaborated for completeness. They concern commutativity of diagrams 
involving Grothendieck's six operations for quasi-coherent complexes.
Readers can skip this section for the first reading. 
 
Let $f\colon X \to Y$ be a map between schemes. For arrows $\rho$ and $\sigma$ in $D(X)$ or $D(Y)$, we write
$
\rho \xLeftrightarrow{A} \sigma
$
to mean $\rho$ and $\sigma$ corresponds to each other under adjunctions indicated by $A$.
We will use $E,F,G,H$ to denote objects in $D(X)$ , and $K,L$ objects in $\rmd(Y)$ or $\rmd(Z)$.
If $\sigma \colon  E \to F$ is an isomorphism in $D(X)$, its  inverse $F \to E$ is 
depicted as $F\xleftarrow[\simeq]{\sigma} E$.
Schemes in the appendix are Noetherian. Unspecified maps in a diagram should be clear from the context.
\subsection{Ad hoc notations}

Most of the material below are lifted from \cite{stacks}.
Let $f\colon X \to Y$ be a map between Noetherian schemes.

\subsubsection{$\otimes_f$}%t
  \label{otimesf}
Given $E,F \in \rmd(X)$, using the adjunction map $f^*f_*\mapsto \id$ for $E$ and $F$,  we obtain $f^*(f_*E\otimes f_*F ) \simeq f^*f_*E \otimes f^*f_*F \to E \otimes F$, it
corresponds under $\lst{f}$ adjunction to a map
\begin{equation}
\tag{$\otimes_f$}
f_*E\otimes f_*F \to f_*(E\otimes F).
\end{equation}

\subsubsection{$\hom_f$} \label{homf}
Consider the composition
\begin{equation} \label{U1}
f_*E \otimes f_*\HOM(E, F) \xrightarrow{\otimes_f} f_*(E \otimes \HOM(E,F)) \xrightarrow{f_*\circ\ev} f_*F.
\end{equation}
By $\tensorhom$ adjunction, we obtain
\begin{equation} \tag{$\hom_f$} 
f_*\HOM(E,F)\to  \HOM(f_*E,f_*F).
\end{equation}

\subsubsection{}
If $G\in \perf(X)$, then $\HOM(G,-) \simeq G^\vee\otimes(-)$ and 
$\HOM(G^\vee,-) \simeq (G^\vee)^\vee\otimes(-)\simeq G \otimes (-)$. 
Therefore $\tensorhom$-adjunction can be rewritten as
\begin{equation}\tag{$\homtensor$}
\HOM(\HOM(G, E), F)\simeq \HOM(E, G \otimes F)
\end{equation}
and 
\begin{equation} \tag{$\tensor$}
\HOM(G^\vee\otimes E, F)\simeq \HOM(E, G \otimes F).
\end{equation}

\subsubsection{$\hom_*^*$}\label{homlst}

For $E\in \rmd(X)$, $K\in \rmd(Y)$, $\lst{f}$ adjunction 
induces a natural isomorphism
\begin{equation}\tag{$\hom_*^*$}
\HOM(K ,f_*E) \to f_*\HOM(f^*K, E).
\end{equation}

Let $L\in D(Y)$, the isomorphism as functors is given by
\[
L \to \HOM(K,f_*E)  \xLeftrightarrow{\tensorhom} L\otimes K \to f_*E
\xLeftrightarrow{\lst{f}} f^*(L\otimes K) \to E,
\] 
and
\[
f^*L\otimes f^*K \simeq f^*(L\otimes K) \to E \xLeftrightarrow{\tensorhom} f^*L\to \HOM(f^*K,E) \xLeftrightarrow{\lst{f}}
L \to f_*\HOM(f^*K,E).
\]

\subsubsection{Projection formula} \label{projfor}

For $E \in \rmd(X)$ and $K \in \rmd(Y)$, the projection formula map
\[ \label{pf}\tag{pf}
f_*E\otimes K \xrightarrow{} f_*(E\otimes f^*K) 
\]
corresponds under $\lst{f}$ adjunction to 
\[
f^*(f_*E\otimes K)\simeq f^*f_*E\otimes f^*K \to E\otimes f^*K,
\] which is the tensor product of the adjunction map
$f^*f_*E\to E$ with $f^*K$.
The map $\text{pf}$ is an isomorphism if $E \in \qcoh(X)$ and  $K\in \qcoh(Y)$
(\cite[\href{https://stacks.math.columbia.edu/tag/08EU}{08EU}]{stacks}).

\subsubsection{$\hom_*^!$} \label{homusk}
Let $f\colon X \to Y$ be a proper map, $E\in \qcoh(X)$, $K\in \qcoh(Y)$. 
Using $\hom_f$ and the adjunction map $f_*f^!\mapsto \id$ we obtain
\begin{equation}\tag{$\hom_*^!$}
f_*\HOM(E, f^!K) \xrightarrow{\hom_f} \HOM(f_*E, f_*f^!K) 
	\xrightarrow{f_*f^!K \to K} \HOM(f_*E, K),
\end{equation}
it corresponds under $\tensorhom$ adjunction
to 
\[
f_*\HOM(E, f^!K)  \otimes f_*E \xrightarrow{\otimes_f} f_*(\HOM(E,f^!K)\otimes E) 
\xrightarrow{f_*\circ \ev} f_*f^!K \xrightarrow{f_*f^*\mapsto\id} K.
\]

When $E\in \coh^{-}(X)$ and $K\in \qcoh^+(Y)$, both $f_*\HOM(E, f^!K)$ and $\HOM(f_*E, f_*f^!K)$ belong to $\qcoh(Y)$
and $\hom_*^!$ is an isomorphism
(\cite[\href{https://stacks.math.columbia.edu/tag/0GEW}{0GEW}]{stacks}).

As a map between functors, for any $L\in \qcoh(Y)$, the isomorphism is described as follows. 
We have \[
L \to f_*\HOM(E, f^!K) \xLeftrightarrow{\lst{f},\tensorhom} f^*L\otimes E \to f^!K \xLeftrightarrow{\usk{f}}
f_*(f^*L\otimes E) \to K, 
\]
and
\[
L \to \HOM(f_*E, K) \xLeftrightarrow{\tensorhom} L\otimes f_*E  \xrightarrow[\simeq]{\mathrm{pf}} f_*(f^*L\otimes E)  \to K,
\] and the two maps on the right hand side are identified.  In general, the description is similar. For $L \in \rmd(Y)$ 
we only have
\[
f^*L\otimes E \to f^!K \xRightarrow{(f_*f^!K \to K) \circ f_*}
f_*(f^*L\otimes E) \to K, 
\] and its composition with the projection formula map
$
L\otimes f_*E  \xrightarrow{\mathrm{pf}} f_*(f^*L\otimes E) 
$ gives $L\otimes f_*E  \to K$.

\subsubsection{${}^*\!\otimes^!$}
Let $K,L \in \qcoh(Y)$.
There is a natural transformation
\begin{equation} \tag{${}^*\!\otimes^!$}
f^*K \otimes f^!L \to f^!(K\otimes L),
\end{equation}
induced by 
\[
f_*(f^*K\otimes f^!L) \xleftarrow[\simeq]{\mathrm{pf}}K\otimes f_*f^!L \xrightarrow{ f_*f^!L \to L} K \otimes L
\]
It is a natural isomorphism if and only if $f_*$ maps $\perf(X)$ to $\perf(Y)$
(\cite[Theorem 5.1,5.4]{Ne}).
In that case, set $L=\cO$, we see that $f^*(-)\otimes \omega_f^\bullet \simeq f^!(-)$, and the
adjunction map $f_*f^!K \to K$ is given by 
\[
f_*f^!K\simeq f_*f^!(K\otimes \cO)\xleftarrow[\simeq]{f_*\circ {}^*\!\otimes^!} f_*(f^*K\otimes \omega_f^\bullet) 
\xleftarrow[\simeq]{\mathrm{pf}}K\otimes f_*\omega_f^\bullet 
\xrightarrow{\id\otimes \tr} K \otimes \cO \simeq K.
\]
Here $\tr$ is the adjunction map $f_*f^!\cO\to \cO$.

\subsection{Identifications results}

\subsubsection{$\hom_f$}
Let $f\colon Z \to X$ be a map of schemes.
We describe the map $\hom_f$ from the functor of points viewpoint.

Let $F \in \rmd(X)$ and $K,L \in \rmd(Z)$, by adjunctions we have
\begin{equation}  \label{B1}
 F \to \HOM(f_*K,f_*L)
 \xLeftrightarrow{\tensorhom}
 F\otimes f_*K \to f_*L, 
\end{equation}
and 
\begin{equation} \label{A1}
F \to  \HOM(f_*K,f_*L) 
\xLeftrightarrow{\tensorhom\ \lst{f}}
f^*(F\otimes f_*K) \simeq f^*F\otimes f^*f_*K \to L.
\end{equation}
The map 
$f_*\HOM(K,L) \xrightarrow{\hom_f} \HOM(f_*K, f_*L)$
corresponds under \eqref{B1} to \eqref{U1}, and under \eqref{A1} to 
\begin{equation}\label{A2}
f^*f_*\HOM(K,L) \otimes f^*f_*K 
\xrightarrow{f^*f_*\mapsto \id}
\HOM(K,L)\otimes K\xrightarrow{\ev} L.
\end{equation}

For any  $E \in \rmd(X)$ and map $\rho\colon  E \to f_*\HOM(K,L)$. From
\[
E  \to f_*\HOM(K,L)\xLeftrightarrow{\lst{f}\ \tensorhom} f^*E\otimes K \to L,
\]
we obtain 
\[ \sigma  \colon f^*E\otimes K \to L\]
as the composition
\begin{equation}   \label{D1}
f^*E\otimes K  
\xrightarrow{f^*\rho \otimes \id } f^*f_*\HOM(K,L)\otimes K 
\xrightarrow{f^*f_* \mapsto \id} \HOM(K,L)\otimes K \xrightarrow{\ev} L.
\end{equation}

Consider the composition 
\begin{equation} \label{Comp}
E \xrightarrow{\rho} f_*\HOM(K,L) \xrightarrow{\hom_f} \HOM(f_*K,f_*L).
\end{equation}
Under \eqref{B1},
it corresponds to 
\[
 E\otimes f_*K \xrightarrow{\rho\otimes \id }
 f_*\HOM(K, L)\otimes f_*K 
 \xrightarrow{\eqref{U1}} f_*L,
\] and can be identified with 
\begin{equation} \label{B2} 
\xymatrix{
 E\otimes f_*K  \ar@[teal][r]^-{\mathrm{pf}}
 	& f_*(f^*E\otimes K)  \ar@[teal][r]^-{f_*\sigma}
		& f_*L,
}
\end{equation}

this can be seen using the commutative diagram
\[
\xymatrix{
E \otimes f_*K \ar[r]^-{\rho\otimes \id}	\ar@[teal][d]^-{\mathrm{pf}}
	& f_*\HOM(K,L)\otimes f_*K \ar[r]^-{\eqref{U1}} \ar[d]^-{\mathrm{pf}}
		& f_*L \\
f_*(f^*E\otimes K) \ar@[teal][r]
	& f_*(f^*f_*\HOM(K,L)\otimes K) \ar@[teal][r]^-{f^*f_*\mapsto \id}
		&  f_*(\HOM(K,L)\otimes K) \ar@[teal][u]^-{f_*\ev}.
}
\]

Under \eqref{A1},  \eqref{Comp} corresponds to
\begin{equation*} %\label{A3}
f^*E\otimes f^*f_*K \xrightarrow{f^*\rho \otimes \id} f^*f_*\HOM(K,L) \otimes f^*f_*K \xrightarrow{\eqref{A2}} L.
\end{equation*}
and can be identified with
\begin{equation}\label{A4}
\xymatrix{
f^*E\otimes f^*f_*K \ar@[teal][rr]^-{f^*f_*K \to K}
	&& f^*E \otimes K \ar@[teal][r]^-{\sigma}
		& L
}
\end{equation}
using the commutative diagram
\[
\xymatrix{
f^*E\otimes f^*f_*K \ar[r]^-{f^*\rho\otimes \id} \ar@[teal][d]^-{f^*f_*K \to K} 
	& f^*f_*\HOM(K,L) \otimes f^*f_*K \ar[d]^-{f^*f_*K \to K}  \ar[r]^-{\eqref{A2}}
		& L\\
f^*E\otimes K 	\ar@[teal][r]^-{f^*\rho\otimes \id}
	& f^*f_*\HOM(K,L)  \otimes K 	\ar@[teal][r]^-{f^*f_*\mapsto \id} 
		&\HOM(K,L)\otimes L \ar@[teal][u]^-{\ev}.
}
\]

\subsubsection{$\hom_f$, $\hom_*^*$ and $\hom_*^!$}

\begin{lemma} \label{dual1}
Let $f\colon Z \to X$ be a map of schemes, and $K,L\in \rmd(Z)$. The composition of 
 \[
 \hom_f\colon
f_*\HOM(K,L) \to \HOM(f_*K, f_*L)
\]
with  
\[\hom_*^* \colon \HOM(f_*K, f_*L) \xrightarrow{\sim} f_*\HOM(f^*f_*K, L)
\] is the same as applying $f_*\HOM(-, L)$ to $f^*f_*K \to K$.

\end{lemma}

\begin{proof}

The composition $F \to \HOM(f_*K,f_*L) \xrightarrow{\sim} f_*\HOM(f^*f_*K,L)$ for any $F \in \rmd(X)$
is determined by
\begin{equation}\label{V1}
F \to \HOM(f_*K,f_*L) \xLeftrightarrow{\eqref{A1}}
f^*F\otimes f^*f_*K \to L \xLeftrightarrow{\tensorhom, f^*\text{-}f_*}
F \to f_*\HOM(f^*f_*K, L).
\end{equation}

For any  $E \in \rmd(X)$ and map $\rho\colon  E \to f_*\HOM(K,L)$, 
the composition 
\[
E \xrightarrow{\rho} f_*\HOM(K,L) \xrightarrow{\hom_f} \HOM(f_*K,f_*L)
\]
 corresponds under \eqref{A1} to \eqref{A4}, 
 then by the right half of \eqref{V1}, we see the map
\[
E \xrightarrow{\rho} f_*\HOM(K,L) \xrightarrow{\hom_f} \HOM(f_*K,f_*L) 
\xrightarrow{\hom_*^*} f_*\HOM(f^*f_*K, L)
\] 
is induced by $f^*f_*K\to K$ under $\lst{f}$ and $\tensorhom$ adjunctions.

\end{proof}

\begin{lemma} \label{dual2}
Let $f\colon Z \to X$ be a map of schemes, $K \in \coh^-(Z), L \in \qcoh^+(Z)$. The composition of
 \[
 \hom_f\colon
f_*\HOM(K,L) \to \HOM(f_*K, f_*L)
\]
 with the inverse of
\[
\hom_*^!\colon 
\HOM(f_*K, f_*L) \xleftarrow{\sim} f_*\HOM(K, f^! f_*L)
\]is the same as
applying $f_*\HOM(K,-)$ to $L \to f^! f_*L$.
\end{lemma}

\begin{proof}

Let $F\in \qcoh(X)$, recall \eqref{B1}
\[
F\to \HOM(f_*K,f_*L)
  \xLeftrightarrow{\tensorhom}
    F\otimes f_*K \xrightarrow{}  f_*L,
\]
and
\begin{equation}\label{E2}
F\to f_*\HOM(K, f^!f_*L) \xLeftrightarrow{f^*\text{-}f_*, \tensorhom} f^*F \otimes K \to f^!f_*L \xLeftrightarrow{f_*\text{-}f^!} f_*(f^*F\otimes K) \to f_*L,
\end{equation}
under the isomorphism $\hom_*^!$ and 
$F\otimes f_*K \xrightarrow{\mathrm{pf}\simeq}f_*(f^*F\otimes K)$
, both map on the left side correspondences to the same map $F\otimes f_*K \to f_*L$.

For any map $\rho \colon E \to f_*\HOM(K,L)$ which corresponds to $ \sigma  \colon f^*E\otimes K \to L$
under
\[
E  \to f_*\HOM(K,L)\xLeftrightarrow{\lst{f}\ \tensorhom} f^*E\otimes K \to L,
\]
the composition of $\rho$ and $\hom_f$ is 
determined under \eqref{B1} by  \eqref{B2}
\[
E \otimes f_*K  
\xrightarrow[\simeq]{\mathrm{pf}}  
f_*(f^*E\otimes K)
\xrightarrow{f_*\sigma}
f_*L.
\]
Consider the right half of \eqref{E2}.  From \eqref{B2},
we get $f^*E\otimes K \to f^!f_*L$ 
as the composition
\[
%\begin{equation}
f^*E\otimes K \xrightarrow{\id\mapsto f^!f_*} f^!f_*(f^*E\otimes K )
\xleftarrow[\simeq]{f^!\circ \mathrm{pf}} f^!(E\otimes f_*K) 
\xrightarrow{f^!\circ \eqref{B2}} f^!f_*L
%\end{equation}
\]
or
\begin{equation}\label{E3}
f^*E\otimes K \xrightarrow{\id\mapsto f^!f_*} f^!f_*(f^*E\otimes K )
\xrightarrow{f^!f_*\sigma} f^!f_*L.
\end{equation}

From the commutative diagram
\[
\xymatrix{
f^*E\otimes K \ar[r]^-{\sigma} \ar[d]^-{\id\mapsto f^!f_*}
	&L \ar[d]^-{\id\mapsto f^!f_*}\\
f^!f_*(f^*E\otimes K) \ar[r]^-{f^!f_*\sigma}
	&f^!f_*L,
}
\]
we see the map \eqref{E3} is the same as
\[
f^*E\otimes K \xrightarrow{\sigma} L \xrightarrow{\id \mapsto f^!f_*}f^!f_*L,
\]
and the lemma follows.

\end{proof}

\subsubsection{perfect objects and adjunctions}

\begin{lemma} \label{lstperf}
Let $f\colon Z \to X$ be a map, $K\in \qcoh(Z)$ and $G\in \perf(X)$, the isomorphism 
\[
\hom_*^*\colon \HOM(G,f_*K) \to f_*\HOM(f^*G, K)
\]
is the composition
\[
 \HOM(G,f_*K) \simeq G^\vee\otimes f_*K \xrightarrow{\mathrm{pf}} f_*(f^*(G^\vee) \otimes K)
 \xrightarrow{(f^*G)^\vee\simeq f^*(G^\vee)} f_*((f^*G)^\vee \otimes K)\simeq f_*\HOM(f^*G,K).
\]
\end{lemma}
\begin{proof}The lemma follows from the construction of the isomorphisms and the fact that
the evaluation and coevaluation maps of $f^*G\in \perf(Z)$ are the pullbacks along $f$ of those of $G$.
\end{proof}

\begin{lemma} \label{dual3}
Let $f\colon Z \to X$ be a proper map, $K \in \qcoh(Z)$ and  $G\in \perf(X), H\in \qcoh(X)$.
The composition of the natural transformations
\[
\begin{array}{cll}
f_*\HOM(\HOM(f^*G, K), f^!H)
&\xrightarrow{ \hom_*^!} &\HOM(f_*\HOM(f^*G, K), H) \\
&  \xrightarrow{(\hom_*^*)^{-1} } & \HOM(\HOM(G, f_*K), H) \\
& \xrightarrow{\homtensor} & \HOM(f_*K, G\otimes H)
\end{array}
\]
is the same as
\[
\begin{array}{cll}
 f_*\HOM(\HOM(f^*G, K), f^!H)  
& \xrightarrow{\homtensor}  &f_*\HOM(K, f^*G\otimes f^! H)\\
& \xrightarrow{{}^*\!\otimes^!} &f_*\HOM(K, f^!(G\otimes H)) \\
 & \xrightarrow{\hom_*^!} & \HOM(f_*K, G\otimes H).
\end{array}
\]
Each transformation $\to$ is induced by certain adjunction whose notation is placed over it.
\end{lemma}

\begin{proof}
Let $E \in \rmd(X)$ and consider a map $\rho\colon E \to f_*\HOM(\HOM(f^*G, K), f^!H)$
which corresponds under $\lst{f}$ and $\tensorhom$ adjunctions to
\begin{equation}\label{J1}
	f^*E\otimes \HOM(f^*G, K) \to  f^!H.
\end{equation} 
We obtain via $\hom_*^!$ in the first sequence of transformations 
\begin{equation}
	E \to \HOM(f_*\HOM(f^*G, K), H),
\end{equation}
it corresponds under $\tensorhom$ adjunction to
\begin{equation}  \label{J2}
	E \otimes  f_*\HOM(f^*G,K) \xrightarrow{\mathrm{pf}} f_*(f^*E \otimes\HOM(f^*G, K)) 
	\xrightarrow{f_*\circ\eqref{J1}} f_*f^!H \xrightarrow{f_*f^!\mapsto\id} H.
\end{equation}
Using $\hom_*^* \colon \HOM(G,f_*K) \simeq f_*\HOM(f^*G, K)$, 
rewrite \eqref{J2} as  
\begin{equation} \label{J3}
	E \otimes  G^\vee \otimes f_*K \simeq E \otimes  \HOM(G,f_*K) \to H,
\end{equation}
then by $\homtensor$ adjunction, we obtained
\begin{equation}\label{J4}
  \begin{array}{rl}
	E \otimes  f_*K  \simeq \cO\otimes E\otimes f_*K 
	&\xrightarrow{\mathrm{coev}\otimes \id} (G\otimes G^\vee)\otimes E\otimes f_*K  
	\simeq G \otimes (E\otimes G^\vee \otimes f_*K) \\
	&\xrightarrow{\id_G\otimes \eqref{J3}} G\otimes H.
  \end{array} 
\end{equation}
It corresponds under $\tensorhom$ adjunction to 
$E \to \HOM(f_*K, G\otimes H)$ obtained from applying the first sequence of transformations
to $\rho$.

For the second sequence of transformations, from \eqref{J1} we obtain
\begin{equation}\label{J5}
\begin{array}{rl}
	f^*E \otimes K 
	\simeq  f^*\cO   \otimes f^*E \otimes K 
	& \xrightarrow{\mathrm{coev}\otimes \id} 
		(f^*G  \otimes (f^*G)^\vee ) \otimes f^*E  \otimes K
		\simeq f^*G  \otimes (f^*E \otimes \HOM(f^*G,K)) \\
	&  \xrightarrow{\id\otimes \eqref{J1}} f^*G\otimes f^!H 
		\xrightarrow{{}^*\otimes^!}  f^!(G\otimes H),
\end{array}
\end{equation}
it corresponds under $\tensorhom,\lst{f}$ adjunctions to 
the map obtained by applying the transformations $\homtensor$ and then ${}^*\!\otimes^!$ 
to $\rho$.
Then via $\hom_*^!$ it corresponds to 
\[
\xymatrix{
  E \otimes f_*K \ar@[teal][r]^-{\mathrm{pf}}
  	& f_*(f^*E \otimes K) \ar@[teal][r]^-{f_*\circ\eqref{J5}}
  		& f_*f^!(G\otimes H)\ar@[teal][r]^-{f_*f^!\mapsto \id}
  			& G\otimes H.
}
\]
This map can be identified with \eqref{J4} using Lemma $\ref{lstperf}$ and the commutative diagram
\[
\xymatrix{
E\otimes f_*K \ar@[teal][rr]^-{\mathrm{pf}} \ar[d]
	&
		&   f_*(f^*E\otimes K)  \ar@[teal][d]\\
G \otimes G^\vee \otimes E\otimes f_*K  \ar[rr]^-{\mathrm{pf}}  \ar[d]^-{\simeq} 
	&
		& f_*(f^*G \otimes f^*G^\vee \otimes f^*E \otimes K) \ar@[teal][d]^-{\simeq}  \\
G\otimes E\otimes f_*\HOM(f^*G, K)  \ar[r]^-{\mathrm{pf}} \ar[ddr]_-{\id_G\otimes \eqref{J2}}
	& G\otimes f_*(f^*E\otimes \HOM(f^*G, K)) \ar[r]^-{\mathrm{pf}} \ar[d]^-{\id\otimes f_*\circ\eqref{J1}}
		& f_*(f^*G\otimes f^*E\otimes \HOM(f^*G, K))  \ar@[teal][d]^-{f_*(\id\otimes \eqref{J1})} \\
	&	G\otimes f_*f^!H \ar[r]^{\mathrm{pf}} \ar[d]^-{f_*f^!H\to H}
		&f_*(f^*G \otimes f^!H) \ar@[teal][d]^-{f_*\circ ^*\otimes^!} \\
	&G\otimes H
		&f_*f^!(G\otimes H) \ar@{}[ul]|{\Box} \ar@[teal][l]^-{f_*f^!\mapsto \id}
}
\]
Note that it follows from the construction of ${}^*\otimes^!$ that the square $\Box$ commutes.

\end{proof}

\subsubsection{pushforward and $\ev$}

\begin{lemma}  \label{dualofpf}
Let $f\colon X \to Y$ be a map, $E\in \qcoh(X)$ and $K\in \qcoh(Y)$. The map
\[
f_*\HOM(E, f^!K) \otimes f_* E 
		\xrightarrow{\otimes_f} f_*(\HOM(E, f^!K) \otimes E )
			\xrightarrow{f_*\circ \ev}f_*f^!K
				\to K
\]
is the same
as \[
f_*\HOM(E, f^!K) \otimes f_* E 
 \xrightarrow{\hom^!_*\otimes \id} \HOM(f_*E, K)  \otimes f_*E 
	\xrightarrow{\ev}   K.
\]
In particular, for $E\in \perf(X)$ and $K=\cO_Y$, 
the isomorphism 
\[
f_*\HOM(E, f^!\cO_Y) \xrightarrow{\hom_*^!}\HOM(f_*E,\cO_Y)
\]
identify the dual of $f_*E$ with $f_*(E^\vee\otimes \omega_f^\bullet)$, and 
the evaluation map
\[
\ev\colon \HOM(f_*E,\cO_Y)\otimes f_*E \to \cO_Y
\] is identified with
\[
f_*\HOM(E, \omega_f^\bullet) \otimes f_*E  \xrightarrow{f_*\ev \otimes_f} f_*\omega_f^\bullet 
=f_*f^!\cO_Y \to  \cO_Y.
\]

\end{lemma}

\begin{proof}
This follows from the construction of $\hom_*^!$ using the first map.
Let $L \in  \rmd(Y)$ and $\rho\colon L \to \HOM(f_*E, K)$, under the correspondence
\[
L \to \HOM(f_*E, K) \xLeftrightarrow{\tensorhom} L\otimes f_*E \to K,
\] the map $L \otimes f_*E \to K$ that corresponds to $\rho$ is
the composition
\[
L\otimes f_*E \xrightarrow{\rho\otimes\id} \HOM(f_*E, K)\otimes f_*E \xrightarrow{\ev} K.
\]

\end{proof}

\subsubsection{$\phi$ and $\psi$}
Let $f\colon Y \to C$ be a proper, flat map, so that $f^*(-)\otimes f^!\cO_Y \xrightarrow{\sim} f^!(-)$.
Let $F \in \perf(Y)$, $K,L\in \qcoh(C)$, 
and $\rho\colon f^*K \to F\otimes f^*L$ a map in $\qcoh(Y)$.

Start with $\rho$, we obtain two maps $\phi$ and $\psi$ as follows.
We have \[
\rho \xLeftrightarrow{\lst{f}}
	K \xrightarrow{}f_*(F\otimes f^*L) \simeq f_*F\otimes L
	 \xLeftrightarrow{\tensor} K\otimes (f_*F)^\vee \xrightarrow{\phi} L, 
\]
and 
\[
\rho \xLeftrightarrow{\tensor}F^\vee\otimes f^*K \to f^*L
	\xRightarrow{ \otimes \omega_f^\bullet} 
  	F^\vee\otimes \omega_f^\bullet \otimes f^*K \to f^*L\otimes \omega_f^\bullet \simeq f^!L
\xLeftrightarrow{\usk{f}} K\otimes f_*(F^\vee \otimes \omega_f^\bullet) \xrightarrow{\psi} L.
\]

\begin{lem} \label{pots}
The map $\phi$ and $\psi$ are the same under $f_*(F^\vee\otimes \omega_f^\bullet) \simeq (f_*F)^\vee $.
\end{lem}
\begin{proof}
Let $\sigma\colon K \to f_*F\otimes L$ be the map corresponding to $\rho$
under $\lst{f}$ adjunction.
We recover $\rho$
as 
\[
f^*K \xrightarrow{f^*\sigma} f^*(f_*F\otimes L)\simeq f^*f_*F\otimes f^*L \xrightarrow{f^*f_*F \to F} F\otimes f^*L.
\]
From $\sigma$, we obtain 
the map $\phi$ as the composition
\[
K\otimes (f_*F)^\vee \simeq \HOM(f_*F, \cO)  \otimes K
 \xrightarrow{\id\otimes\sigma} \HOM(f_*F,\cO)\otimes f_*F \otimes L \xrightarrow{\ev\otimes \id_L} \cO\otimes L \simeq L.
\]
By Lemma \ref{dualofpf}, it is the same as
\[
 f_*\HOM(F, \omega_f^\bullet)  \otimes K
 \xrightarrow{\id\otimes\sigma} 
 f_*\HOM(F, \omega_f^\bullet) \otimes f_*F \otimes L \
 \xrightarrow{f_*\ev\circ \otimes_f \ \otimes \id_L} f_*\omega_f^\bullet \otimes L \
 \xrightarrow{\tr\otimes\id} \cO\otimes L\simeq L.
\] 
Identify it with the map $\psi$ using the commutative diagram
\[
\xymatrix{
f_*\HOM(F, \omega_f^\bullet) \otimes K \ar[d]^-{\id\otimes \sigma} \ar@[teal][r]^-{\mathrm{pf}}
	& f_*(\HOM(F, \omega_f^\bullet) \otimes f^*K) 
		\ar@/^5pc/@< 12pt>@[teal][dd]^-{f_*\circ (\id\otimes \rho)}
		\ar[d]^-{f_*\circ(\id\otimes f^*\!\sigma)}
	  \\
 f_*\HOM(F, \omega_f^\bullet)\otimes f_*F\otimes L \ar[r]^-{\mathrm{pf}} \ar[d]^-{\otimes_f\ \otimes \id_L}
 	&f_*(\HOM(F, \omega_f^\bullet)\otimes f^*(f_*F\otimes L ))\ar[d]^-{f^*f_*F\to F} 
		\\
f_*(\HOM(F,\omega_f^\bullet)\otimes F)\otimes L  \ar[r]^-{\mathrm{pf}} \ar[d]^-{f_*\circ\ev\otimes \id_L}
	& f_*(\HOM(F, \omega_f^\bullet)\otimes F\otimes f^*L) \ar@[teal][d]^-{f_*\circ (\ev\otimes \id)}
		\\
 f_*\omega_f^\bullet\otimes L  \ar[r]^-{\mathrm{pf}}  \ar[d]^-{\tr\otimes \id_L}
 	& f_*(\omega_f^\bullet\otimes f^*L)	\ar@[teal][d]_-{\simeq}^-{^*\!\otimes^!}\\
 L	&	f_*f^!L \ar@[teal][l]^-{f_*f^!L \to L}
}
\]
Note that the composition of the colored arrows gives $\psi$.
\end{proof}

\end{appendix}

\providecommand{\bysame}{\leavevmode\hbox to3em{\hrulefill}\thinspace}
\providecommand{\MR}{\relax\ifhmode\unskip\space\fi MR }
\providecommand{\MRhref}[2]{%
 \href{http://www.ams.org/mathscinet-getitem?mr=#1}{#2}}
\providecommand{\href}[2]{#2}

\end{document}